\documentclass[12pt,reqno]{amsart}

\usepackage[LGR,T1]{fontenc}
\usepackage[utf8]{inputenc}

\title[Time reversal of Markov processes with jumps]{Time reversal of Markov processes with jumps under a finite entropy condition}
\date{October 2021, revised version}

\author[Conforti]{Giovanni Conforti}
\address{G.~Conforti. Département de Mathématiques Appliquées, \'Ecole Polytechnique, Palaiseau. France}
\email{giovanni.conforti@polytechnique.edu}

\author[Léonard]{Christian Léonard}
\address{C.~Léonard. Modal’X, UPL, Université Paris Nanterre. France}
\email{christian.leonard@math.cnrs.fr.fr}
\thanks{This research is partially  granted by the projects   SPOT (ANR-20-CE40-0014) and  Labex MME-DII (ANR-11-LBX-0023)}

\usepackage[english]{babel}

\usepackage{amssymb, amsmath, amsfonts, latexsym, enumerate}
\usepackage[usenames,dvipsnames]{pstricks}
\usepackage{epsfig}
\RequirePackage[colorlinks,linkcolor=blue,citecolor=blue,urlcolor=blue]{hyperref}

\newtheorem{theorem}[equation]{Theorem}
\newtheorem{lemma}[equation]{Lemma}
\newtheorem{proposition}[equation]{Proposition}
\newtheorem{corollary}[equation]{Corollary}

\newtheorem{hypotheses}[equation]{Hypotheses}

\theoremstyle{remark}
\newtheorem{remark}[equation]{Remark}
\newtheorem{remarks}[equation]{Remarks}

\numberwithin{equation}{section}

\newcommand{\RR}{\mathbb{R}}
\newcommand{\Rn}{\mathbb{R}^n}

\newcommand{\1}{\mathbf{1}}

\newcommand{\ttimes}{\!\times\!}

\newcommand\pf{_{\#}}

\newcommand{\Leb}{\mathrm{Leb}}

\renewcommand{\ae}{\textrm{-}\mathrm{a.e.}}
\newcommand{\Id}{\mathrm{Id}}
\newcommand{\scal}{\!\cdot\!}
\newcommand*{\cchi}{\raisebox{0.35ex}{\( \chi \)}}

\DeclareMathOperator{\dom}{dom}

\DeclareMathOperator{\supp}{supp}

\DeclareMathOperator{\MP}{MP}

\newcommand{\boulette}[1]{$\bullet$\ Proof of #1.}
\newcommand{\Boulette}[1]{\par\medskip\noindent $\bullet$\ Proof of #1.}
\newcommand{\sbt}{\,\begin{picture}(-1,1)(-1,-3)\circle*{3}\end{picture}\ }

\newcommand\Lim[1]{\lim_{#1\rightarrow\infty}}

\newcommand\Limh{\lim_{h\to 0^+}}

\newcommand\II[2]{\int_{[#1,#2]}}

\newcommand{\cadlag}{c\`adl\`ag}

\newcommand{\ii}{[0,T]}
\newcommand{\io}{[t_o,T]}
\newcommand\XX{ \mathcal{X}}

\newcommand\iX{{\ii\ttimes\XX}}

\newcommand{\ZZ}{\Rn}
\newcommand\ZZZ{(\ZZ)^2}

\newcommand{\iZ}{\ii\times\ZZ}
\newcommand{\ioZ}{\io\times\ZZ}
\newcommand{\Zs}{\ZZ_\ast}
\newcommand{\iZZ}{\iZ\times\ZZ}

\newcommand{\iZZs}{\iZ\times\Zs}
\newcommand{\ioZZs}{\ioZ\times\Zs}

\newcommand\OO{\Omega}

\newcommand{\UU}{ \mathcal{U}}

\newcommand{\UQ}{ \UU^Q}

\newcommand\PX{\mathrm{P}(\XX)}

\newcommand\MX{\mathrm{M}(\XX)}

\newcommand\PZ{\mathrm{P}(\ZZ)}

\newcommand\MZ{\mathrm{M}(\ZZ)}
\newcommand\PO{\mathrm{P}(\OO)}
\newcommand\MO{\mathrm{M}(\OO)}
\newcommand{\CcZ}{C^2_c(\ZZ)}

\newcommand\Iii{\int_{\ii}}

\newcommand\IZ{\int_{\ZZ}}
\newcommand\IZZ{\int_{\ZZZ}}
\newcommand\IiZ{\int_{\iZ}}
\newcommand\IiZZ{\int_{\ii\times\ZZZ}}

\newcommand{\IZs}{\int _{ \Zs}}
\newcommand{\IZZs}{\int _{\ZZ\times \Zs}}
\newcommand{\IiZZs}{\int _{\ii\times\ZZ\times \Zs}}

\newcommand\Xb{\overline{X}}

\newcommand\Qb{\overline{Q}}

\newcommand{\LL}{ \mathcal{L}}
\newcommand{\LLf}{\overrightarrow{ \mathcal{L}}}
\newcommand{\LLb}{\overleftarrow{ \mathcal{L}}}

\newcommand{\Gf}{\overrightarrow{ \mathcal{G}}}
\newcommand{\Gb}{\overleftarrow{ \mathcal{G}}}
\newcommand{\Gaf}{\overrightarrow{\Gamma}}
\newcommand{\Gab}{\overleftarrow{\Gamma}}
\newcommand{\Jf}{\overrightarrow{J}}
\newcommand{\Jb}{\overleftarrow{J}}

\newcommand{\If}[2]{[#2\overrightarrow{J} ]_{ #1}}
\newcommand{\Kf}{\overrightarrow{K}}
\newcommand{\Kb}{\overleftarrow{K}}

\newcommand{\Bf}{\overrightarrow{b}}
\newcommand{\Bb}{\overleftarrow{b}}
\newcommand{\Btf}{\overrightarrow{b} ^{ \delta}}
\newcommand{\Btb}{\overleftarrow{b} ^{ \delta}}

\newcommand{\pp}{ \mathsf{p}}
\newcommand{\qq}{ \mathsf{q}}
\newcommand{\rr}{ \mathsf{r}}

\newcommand{\hh}{ \mathsf{h}}

\newcommand{\pb}{\bar{\pp}}
\newcommand{\qb}{\bar{\qq}}

\newcommand{\flo}[1]{\lfloor #1\rfloor ^{ \delta}}
\newcommand{\floo}[2]{\lfloor #1-#2\rfloor ^{ \delta}}
\newcommand{\xd}{\lfloor \xi \rfloor^ \delta}
\newcommand{\FE}{ \mathrm{FE}}
\newcommand{\IZA}{\int _{ \ZZ\times \mathcal{A}}}
\newcommand{\faa}{\hat\phi_\alpha}
\newcommand{\fa}{\phi_\alpha}

\begin{document}

\begin{abstract} 
Motivated by entropic optimal transport, time reversal of Markov jump processes in $\ZZ$  is investigated. 
Relying on an abstract  integration by parts formula  for the carré du champ of a  Markov process recently  obtained in \cite{CCGL20}, and using an entropic improvement strategy discovered by Föllmer \cite{Foe85b,Foe86},   we compute the semimartingale characteristics of the time reversed   process   for a wide class of jump  processes in $\ZZ$  with possibly unbounded variation sample paths and singular  intensities of jump. 
\end{abstract}

\maketitle 
\tableofcontents

\section{Introduction}

The time-reversed   of a Markov process   remains a Markov process. Consequently, the problem of finding its Markov generator  arises. The answer to this problem is given by a   time reversal formula.

To our knowledge,  general results in terms of semimartingale characteristics  of time-reversed Markov processes with jumps in a continuous time setting  are not available in the literature. Other types of results are known,  for instance  Jacod and Protter identify in \cite{JP88} a  large class of semimartingales built upon Lévy processes which remain semimartingales once time-reversed. This is a nontrivial result because Walsh showed in \cite{W82} that time-reversing a semimartingale might not give a semimartingale anymore. The present article exhibits a large class of Markov semimartingales  with jumps whose bounded variation parts are absolutely continuous with a  possibly singular derivative   (the intensity of jump) and whose time-reversed  are still semimartingales with an absolutely continuous part.

 Of course, the intuition for the expression \eqref{eq-65pb} below of the  jump intensities of the time reversal of a process with jumps is strong. It appears at the very beginning of the story in Kolmogorov's celebrated article \cite[Eq.\,(7)]{Kol36}, and any  physicist writes it without hesitating. Nevertheless, a complete proof for rather general processes with possibly singular intensities of jumps was not done. The present article provides such results.

\subsubsection*{Entropic optimal transport}

Besides being an interesting topic in its own right,  last years have seen a renewed interest in  time reversal  because of its applications to entropic optimal transport (Schrödinger problem) and functional inequalities.  We refer to the seminal article \cite{Zam86} by Zambrini, the survey paper \cite{Leo12e} by the second author for more about this active research field.

Time reversal is invoked to study entropic optimal transport in the recent articles  by  Chen, Georgiou and Pavon \cite{CGP14}, Gentil, Léonard and Ripani \cite{GLR17a}, Conforti \cite{Co18}, Backhoff, Conforti, Gentil and Léonard \cite{BCGL19}. 
Regarding functional inequalities,  the   logarithmic Sobolev  and  HWI inequalities, and  the Bakry-\'Emery criterion are recovered  using   time reversal by  Fontbona and  Jourdain \cite{FJ16}, Gentil, Léonard and Ripani \cite{Leo12d,GLR17b} and   Karatzas, Schachermayer and Tschiderer \cite{KST18}. 

As exposed in  the introductory section of \cite[p.\,3]{CCGL20} about entropic and deterministic optimal transports, in the setting of diffusion processes  time reversal allows for an expression of the entropy of a diffusion measure as the difference of a kinetic action on  the Wasserstein space of probability measure, and a potential action whose integrand is minus some Fisher information. It is a work in progress to derive from this classical mechanical picture a Newton equation for the entropic interpolations in the spirit of the articles \cite{vR11} by von Renesse and \cite{Co18} by the first author.  

All these contributions take place in a diffusion setting.  
Their  analogues in presence  of jumps remain to be explored. This article   is   a preliminary step to obtain a similar Lagrangian action representation of the relative entropy of a path measure in the setting of  jump processes. 

To our knowledge the only articles where  jump processes are investigated in the context of the Schrödinger problem are \cite{PZ04,PZ05} by Privault and Zambrini. Although not referring directly to any time reversal formula, time reversal plays a crucial role in these papers.

\subsubsection*{Backward Fokker-Planck and Hamilton-Jacobi-Bellman equations}

In a different context,   Izydorczyk, Oudjane, Russo and  Tessitore \cite{IOR21,IORT21} recently used time reversal of diffusion processes  to prove the well-posedness of some backward Fokker-Planck equations, and to design  efficient algorithms  solving some Hamilton-Jacobi-Bellman equations with terminal conditions. It is natural to ask for  analogous results when replacing the second order elliptic term in these PDEs  by a nonlocal Markov generator attached to a jump process. The present article provides the necessary time reversal formula.

\subsection*{Main results of the article}
Our main results are Theorems \ref{res-43} and \ref{res-37}.
\\
Let us briefly present the  content of Theorem  \ref{res-43} skipping its  detailed hypotheses, in the simple case where the sample paths have bounded variation. Its proof is partly based on the abstract time reversal formula of Theorem \ref{res-37}. Consider a Markov process with  generator  defined for any function $u$ in $ C_c^1(\ZZ)$, by
\begin{align*}
 \Bf(t,x)\scal\nabla u(x)+\IZs [u(y)-u(x)]\, \Jf _{ t,x}(dy),
\quad (t,x)\in\iZ,
\end{align*}
where $\Bf$ is a vector field and the jump kernel $\Jf$ satisfies $\IZs (|y-x|\wedge 1)\,\Jf _{ t,x}(dy)< \infty$ for any $t,x$ to assure that the  integral  in the expression of the generator is well-defined. This also implies that the sample paths  have bounded variation. 
\\
Then, under some hypotheses, its time-reversed process admits the Markov generator defined by
\begin{align*}
\ \Bb(T-t,x)\scal\nabla u(x)+\IZs [u(y)-u(x)]\, \Jb _{T- t,x}(dy),
\quad (t,x)\in\iZ,
\end{align*}
with 
\[\Bb _{t}=-\Bf_t\]   and the backward  jump kernel $\Jb$  is the unique solution of  the flux equation
\begin{align}\label{eq-65pb}
\pp_t(dy) \Jb _{ t,y}(dx)=\pp_t(dx)\Jf _{ t,x}(dy),
\end{align}
where the known coefficients are the forward jump kernel $\Jf$ and the  marginal $\pp_t$: the law of the forward process at time $t$.
\\
Formula \eqref{eq-65pb} remains valid when the sample paths  are only supposed to have bounded quadratic variation. This  happens when     
\begin{align*}
\IZs (|y-x|^2\wedge 1)\,\Jf _{ t,x}(dy)< \infty,
\quad (t,x)\in\iZ.
\end{align*}
In this wider setting, the expression of the generator requires some truncation technicalities from which we stay apart during this introduction, see \eqref{eq-60} below.

Typically,  Theorem \ref{res-43}  is proved assuming that $\Bf$ is a regular vector field. But no regularity is required for the jump kernel except some ``entropic integrability'' (see  Corollary \ref{res-44}) allowing \emph{locally unbounded intensities of jumps}. It also states that if the time reversal formula holds for some reference Markov measure $R$, then it also holds for any Markov measure $P$ such that its relative entropy 
\[
H(P|R):=E_P \log( dP/dR)< \infty
\] 
with respect to $R$ is finite. \emph{This is precisely what is needed for entropic optimal transport} where no a priori regularity is known except this finite entropy estimate.

These results are consequences of
\begin{enumerate}
\item
 the abstract time reversal formula of Theorem \ref{res-37} which permits to obtain at
 \item
  Theorem \ref{res-36} a  first time reversal formula under the hypothesis that the jump kernel is regular, and then to extend it to singular kernels using  
 \item
 the entropic improvement   Lemma \ref{res-18}.
\end{enumerate}

\subsection*{Literature}
As already alluded to, the literature on this topic is sparse.  
The time reversal formula \eqref{eq-65pb} is similar to the one obtained for random walks on graphs by Cattiaux, Gentil and the authors  in \cite{CCGL20} which is a paper  whose main concern is to extend  already known   time reversal formulas for diffusion processes \cite{Foe85b,Foe86,Par86,HP86,MNS89} to a wide class of diffusion processes with singular drifts, again with entropic optimal transport in mind. For a little more about the literature on time reversal of Markov processes, one may have a look at the introduction of \cite{CCGL20}.

As in \cite{CCGL20}, the idea of the entropic improvement at Section \ref{sec-entropy}, leading us in the present framework to hypotheses allowing singular  intensities, is fully credited to Föllmer \cite{Foe85b,Foe86}.

\subsection*{Outline of the article}
Our  results  rely on an integration by parts (IbP) formula  for the carré du champ of a Markov process which was proved by  Cattiaux, Gentil and the authors in \cite{CCGL20}. This IbP formula which is recalled at Section \ref{sec-IbP} serves us to prove at Section \ref{sec-jp-abs} the abstract time reversal formula of Theorem \ref{res-37}.  At Section \ref{sec-regular-kernel}, we apply this abstract result to prove the time reversal formula of Theorem \ref{res-36} in the the case where the jump kernel is regular. This result is extended at Section \ref{sec-entropy} where the proof of Theorem \ref{res-43} is completed by means of an entropic improvement.
\\
The main time reversal formula of this article and the  sets of hypotheses are stated at next Section \ref{sec-hyp}. Examples are displayed at Section \ref{sec-examples} to illustrate the  generality of our assumptions.

\subsection*{Notation}

The set of all probability measures on a measurable set $A$ is denoted by $ \mathrm{P}(A)$ and the set of all nonnegative $ \sigma$-finite measures on $A$ is  $ \mathrm{M}(A).$
The push-forward of a measure $ \qq\in  \mathrm{M}(A)$ by the measurable map $f:A\to B$ is   $f\pf\qq(\sbt):=\qq(f\in \sbt)\in \mathrm{M}(B).$ 

\subsubsection*{Relative entropy}

The relative entropy of $\pp\in \mathrm{P}(A)$ with respect to the reference  measure $\rr\in \mathrm{M}(A)$ is 
\begin{align*}
H(\pp|\rr):=\int_A \log(d\pp/d\rr)\,d\pp \in( - \infty, \infty]
\end{align*}
if $\pp$ is absolutely continuous with respect to $\rr$ ($\pp\ll\rr$) and $\int_A \log_-(d\pp/d\rr)\,d\pp< \infty$, { and $H(\pp|\rr)=+ \infty$ otherwise}. If $\rr\in \mathrm{P}(A)$ is a probability measure, then $H(\pp|\rr) \in[0, \infty].$ See \cite[App.\,B]{CCGL20} for details.

\subsubsection*{Path measures}
The configuration space is a Polish space $\XX$ equipped with its Borel $ \sigma$-field.
The path space is the set $\OO:=D(\ii,\XX)$ of all $\XX$-valued \cadlag\  trajectories on the time index set $\ii,$ and the canonical process  $(X_t) _{ 0\le t\le T}$ is defined by $X_t( \omega)= \omega_t$ for any $0\le t\le T$ and any path $ \omega=( \omega_s) _{ 0\le s\le T}\in\OO.$ It is equipped  with the canonical $ \sigma$-field $ \sigma(X _{ \ii})$ and the the canonical filtration $ \Big(\sigma( X _{ [0,t]}); 0\le t\le T\Big)$ where for  any subset $ \mathcal{T}\subset \ii$, $X _{ \mathcal{T}}:=(X_t, t\in \mathcal{T})$ and  $ \sigma(X _{ \mathcal{T}})$ is the $ \sigma$-field generated by the collection of maps $(X_t, t\in \mathcal{T})$.  
\\
 We call any positive measure  $Q\in\MO$ on $\OO$ a path measure. For any $\mathcal{T}\subset\ii,$ we denote  $Q_\mathcal{T}=(X_\mathcal{T})\pf Q.$ In particular, for any
 $0\le r\le s\le T,$ $X_{[r,s]}=(X_t)_{r\le t\le s}$, $Q_{[r,s]}=(X_{[r,s]})\pf Q$, and $Q_t=(X_t)\pf Q\in\MX$ denotes the law of the position $X_t$ at time $t$. If $Q\in\PO$ is a probability measure, then $Q_t\in\PX$.
\\
The time-space canonical process is
$$
\Xb_t:=(t,X_t)\in\iX,
$$
and for  any function $u:\iX\to\RR$,  we denote  $u(\Xb): (t, \omega)\mapsto u(t, \omega_t).
$
We also denote
\begin{align*}
\Qb  (dtd \omega)&:=dtQ (d \omega),\qquad dtd \omega\subset\ii\times \OO,\\
\qb (dtd x)&:=dtQ_t (d x),\qquad dtd x\subset\iX.
\end{align*}

\section{Main result. Hypotheses}
\label{sec-hyp}

\subsection*{Basic definitions}

Before describing the hypotheses and stating our  main result at Theorem \ref{res-43}, we recall some basic definitions.

\subsubsection*{Conditionable path measure}

  A path measure $Q$ such that  $Q_t$ is  $ \sigma$-finite for all $t$ is called a conditionable path measure. This notion is necessary to define properly the conditional expectations $E_Q(\sbt\mid X _{t }),$ $E_Q(\sbt\mid X _{[0,t] })$ and $E_Q(\sbt\mid X _{[t,T] }),$ for any $t$. If $Q$ has a finite mass, then it is automatically conditionable.

\subsubsection*{Extended forward generator}

 Let $Q$ be a conditionable measure. A measurable function $u$ on $\iX$ is said to be in
the domain of the extended forward generator of $Q$ if there exists a real-valued 
process $\LLf^Q u(t,X _{ [0,t]})$  which is adapted with respect to the forward filtration such that
 $\Iii|\LLf^Q u(t,X _{ [0,t]})|\,dt<\infty,$ $Q\ae$ and the process $$M^u_t:=u(\Xb_t)-u(\Xb_0)-\II0t \LLf^Q u(s,X _{ [0,s]})\,ds,\quad 0\le t\le T,$$ is a local $Q$-martingale.
We say that $\LLf^Q $ is  the extended forward generator of $Q.$ Its domain  is denoted by $\dom
\LLf^Q .$

\subsubsection*{Reversing time}

Let $Q\in\MO$ be any path measure. Its time reversal is 
\begin{align*}
Q^*:=(X^*)\pf Q\in\MO,
\end{align*}
where 
\begin{align*}
\left\{ \begin{array}{ll}
X^*_t:= \Limh X _{ T-t+h},\quad & 0\le t< T,\\
 X^*_T:=X_0,& t=T,
\end{array}\right.
\end{align*}
is the reversed canonical process. We assume that $Q$ is such that $Q(X _{ T^-}\neq X_T)=0,$ i.e.\ its sample paths are left-continuous at $t=T.$ This implies that the time reversal mapping $X^*$ is $Q\ae$\,one-one on $\OO.$ 
 We introduce the backward extended generator 
 \begin{align}\label{eq-F05}
\LLb^Qu(t, X _{ [t,T]})&:=\LLf ^{ Q^*} u^*(t^*,X^* _{ [0,t^*]}),
\end{align}
where $ u^*(t^*, \omega^* _{ [0,t^*]}):=u(t, \omega _{ [t,T]}),$  with $t^*:=(T-t)^+$, $ \omega^*(t):= \omega(t^*)$, and $ \LLf ^{ Q^*}$  stands  for the standard (forward) generator  of $Q^*$.

  \subsubsection*{Markov measure}

  A path measure $Q\in\MO$ is said to be Markov if it is conditionable and for any $0\le t\le T,$ $Q(X _{ [t,T]}\in\sbt\mid X _{ [0,t]})=Q(X _{ [t,T]}\in\sbt\mid X _t).$  It is known that $ Q^*$
 is also Markov and its extended generators at time $t$ only depend of the present position $X_t$. Therefore it is possible to consider the sum and difference of the forward and backward generators: they remain functions of the present position.

%\subsubsection*{Osmotic generator} In restriction to  $\dom\LLf^Q\cap \dom\LLb^Q,$ we define the \emph{osmotic} extended generator of $Q$ by \begin{align*} \LLo Q:= (\LLf^Q+\LLb^Q)/2. \end{align*} It plays an important role in this article. 

\subsection*{Jump process on $\ZZ$}
Let us recall basic notions about jump processes on $\ZZ.$

\subsubsection*{Generator}

The  generator of a general Markov jump process  on $\ZZ$ without diffusion is
\begin{align}\label{eq-60}
\LL_tu(x)
	= b^ \delta (t,x)\scal\nabla u(x)+\IZs [u(x+\xi)-u(x)-   \nabla u(x)\scal \xd]\, K _{ t,x}(d\xi),
\end{align}
where   $K$ is a measurable field of  nonnegative measures on the jump set  $\Zs:=\ZZ\setminus \left\{0\right\} $ such that 
\begin{align}\label{eq-872}
\IZs (|\xi|^2\wedge 1)\, K _{ t,x}(d\xi)< \infty,
\quad \forall (t,x)\in\iZ,
\end{align}
$b^ \delta$ is a locally bounded measurable vector field,
\begin{align*}
\left\{
\begin{array}{ll}
\lfloor\xi\rfloor ^{0}:=0, & \textrm{when } \delta=0;\\
\xd:=\1 _{ \left\{|\xi|\le  \delta\right\} }\,\xi,\qquad 
& \textrm{when } \delta>0,
\end{array}
\right.
\end{align*}
and $u$ belongs to the class 
$C_c^2(\ZZ)$ 
 of twice continuously differentiable functions with a compact support in $\ZZ.$

%$\xd:=[\1 _{ \left\{|\xi|\le  \delta/2\right\} }+\1 _{ \left\{ \delta/2<|\xi|\le  \delta\right\} }(2-2|\xi|/ \delta)]\, \xi.$ As a function of $\xi$ with $0<  \delta\le 1$ fixed, $\xd$  is a continuous truncation function   which is equal to the identity on a small ball (with radius $ \delta/2$) around the origin and zero outside a larger small ball (with radius $ \delta$). 

The truncation $\xd$ with  $ \delta>0$ appears for  the integral in the definition of $\LL_t u$ to be well defined  under the assumption \eqref{eq-872}. %It is useless  as soon as the stronger assumption \eqref{eq-31} below is satisfied; in this case, one takes $ \delta=0$.

\subsubsection*{Variation of the sample paths}
Under \eqref{eq-872} the sample paths have almost surely bounded quadratic variation. 
If \eqref{eq-872} is reinforced by
\begin{align}\label{eq-871}
\IZs (|\xi|\wedge 1)\, K _{ t,x}(d\xi)< \infty,
\quad \forall (t,x)\in\iZ,
\end{align}
then the sample paths have bounded variation. In this case, one chooses $ \delta=0$ in the expression of $\LL,$ leading to the meaningful formula
\begin{align}\label{eq-60b}
\LL_tu(x)
	= b(t,x)\scal\nabla u(x)+\IZs [u(x+\xi)-u(x)]\, K _{ t,x}(d\xi),
\end{align}
with the simplified notation $$b:=b ^{ \delta=0}.$$ When \eqref{eq-871} fails, the integral on the right hand side of \eqref{eq-60b} is undefined and we have to take $ \delta>0,$ for instance $ \delta=1.$ Under \eqref{eq-871}, we see that $b=b^ \delta-\IZs   \xd \, K _{ t,x}(d\xi)$ for any $ \delta\ge 0$, showing that $b^ \delta$ is an artefact which is only necessary when \eqref{eq-871} fails.

\subsubsection*{Martingale problem}

We say that the path measure $Q\in\MO$ solves the martingale problem
\begin{align*}
Q\in\MP _{ \delta}(\qq_0, \Btf, K)
\end{align*}
when $Q_0=\qq_0,$  and for almost every $t$ and  any $u\in\CcZ,$ $\dom\LLf_t^Q$ contains $\CcZ$  and $\LLf^Q_tu=\LL_t u$, see \eqref{eq-60}.
When $K$ satisfies \eqref{eq-871}, we choose 
$ \delta=0$.
\\
When the time interval is $\io $ with $0\le t_o< T,$ we write
\begin{align*}
Q _{ [t_o,T]}\in\MP _{t_o, \delta}(\qq _{ t_o}, \Btf, K),
\end{align*}
meaning that $Q _{ t_o}=\qq _{ t_o}$ and $\LLf_t^Qu=\LL_t u$ for almost every $t\in\io$ and $u\in\CcZ.$
\\
For simplicity, we write 
$Q\in\MP _{t_o, \delta}( \Btf, K)$ instead of $Q\in\MP _{t_o, \delta}(Q _{ t_o}, \Btf, K).$

\subsubsection*{Kernel in terms of jumps or positions}

The  jump  kernel $K _{ t,x}(d\xi)$ which is expressed in terms of the jump $\xi$  can equivalently be expressed in terms of the position $y$ after the jump, leading to
\begin{align*}
\LL_tu(x)
	= b^ \delta (t,x)\scal\nabla u(x)+\IZs [u(y)-u(x)-   \nabla u(x)\scal \floo yx]\, J _{ t,x}(dy),
\end{align*}
where the kernel  $J$ is defined  for any $t$ and $x$ by
\begin{align*}
\IZs f(x+\xi)\, K _{ t,x}(d\xi)=\int _{ \ZZ\setminus \{x\}} f(y)\, J _{ t,x}(dy),
\quad \forall f\in C_c(\ZZ\setminus \{x\}).
\end{align*}

Our  main results: Theorem \ref{res-43} and its Corollary \ref{res-44},  require some regularity and/or  integrability properties of the drift field, the jump kernel and also the time marginals of the Markov measure. 

\subsection*{Hypotheses}

Their hypotheses are  built with some properties which are picked up from a list    labeled from  \eqref{eq-121}  to \eqref{eq-89} that  is postponed after the statements of these results, for a better readability because this list is rather long.

\begin{hypotheses}\ \label{ass-08}
\begin{enumerate}[(a)]
\item \emph{General hypotheses.}\  
		The \emph{drift} field $\Bf ^{\delta}$ is locally bounded: \eqref{eq-121}, 
	the \emph{jump kernel} $K$ satisfies \eqref{eq-32} and 
	the \emph{path measure}   satisfies \eqref{eq-94}.

\item One of the following  hypotheses  is fulfilled. 
	\begin{enumerate}[1-]
	\item\emph{Large jumps 1.}\ 
	 Assumption \eqref{eq-76}.
	 \item\emph{Large jumps 2.}\ 
	 Assumption \eqref{eq-110}.
	\end{enumerate}

\item One of the following sets of hypotheses  is fulfilled.
	\begin{enumerate}[1-]
	\item \emph{Bounded variations.}\ 
	Assumptions: \eqref{eq-31a}, \eqref{eq-771}.
	\item \emph{Unbounded variations.}\ 
	If \eqref{eq-97} holds, assume: \eqref{eq-772}, \eqref{eq-100}-\eqref{eq-119}-\eqref{eq-100a}-\eqref{eq-100b}-\eqref{eq-100c}-\eqref{eq-100d} for the small jumps, and \eqref{eq-88} for the marginal flow.
	\end{enumerate}

\item
One of the following hypotheses is fulfilled. 
	\begin{enumerate}[1-]
	\item
	\emph{Close to reversibility.} Assumption \eqref{eq-81}.
	\item
	\emph{Marginal flow.}\ 
	Assumption \eqref{eq-89}.
	 \item
	 For any $u\in \CcZ,$ $\LLb u$ is bounded.
	 \end{enumerate}
\end{enumerate}
\end{hypotheses}

There are twelve chains of hypotheses: $\{a\}\times \{b1, b2\}\times\{c1,c2\}\times \{d1, d2, d3\}$.

\subsection*{Time reversal formula}

Take a reference Markov measure $R\in\MO$ solving
$
R\in \MP _{  \delta}(\Bf ^{ R, \delta},\Jf^R)
$
and satisfying the Hypotheses \ref{ass-08}. Then, consider another Markov probability measure $P\in\PO$ with a finite entropy with respect to $R$: 
$H(P _{ \io}|R _{ \io}) < \infty,$
for some $0\le t_o< T.$
Next results give time reversal formulas for the restriction $P _{ \io}$ of $P$ to the $ \sigma$-field $ \sigma(X _{ \io})$. 

We   require in addition that $R$ is the unique solution to its own martingale problem   in the following sense
\begin{align}\label{eq-115}
[R'\ll R\ \textrm{and}\ R'\in\MP _{ \delta}(R_0,\Bf ^{ R, \delta},\Jf^R)] \implies R'=R.
\end{align}
For instance, it is known that \eqref{eq-115} is satisfied when $R$ is the law of the unique strong solution of an SDE, see \cite{Jac79}.

\begin{theorem}\label{res-43}
Assume that $R\in\MO$ solves
\begin{align*}
R\in \MP _{  \delta}(\Bf ^{ R, \delta},\Jf^R)
\end{align*}
and satisfies the Hypotheses \ref{ass-08} and \eqref{eq-115}. Suppose also  that for some $0\le t_o< T,$ $P\in\PO$ has a finite entropy with respect to $R$ on the time interval $[t_o,T]$:
\begin{align}\label{eq-116}
H(P _{ \io}|R _{ \io}) < \infty.
\end{align}
Then, 
\begin{align*}
P _{ \io}\in\MP _{ t_o,\delta}(\Bf ^{ P, \delta},\Jf^P)
\end{align*}
for some jump kernel $\Jf^P$, with 
\begin{align}\label{eq-114a}
\begin{split}
\Bf ^{ P, \delta}_t(x)=\Bf ^{ R, \delta}_t(x) + \IZs \floo yx\ (\Jf^P _{ t,x}-\Jf^R _{ t,x})&(dy),\\
 &(t,x)\in\ioZ, \ \pb\ae
\end{split}
\end{align}
Moreover, for almost every $t\in\io$,   any $u\in\CcZ$ is in $ \dom \LLb^P_t$ and 
\begin{align}\label{eq-677}
\begin{split}
\LLb^P_t u(x)= \Bb ^{ P, \delta}_t(x)\cdot \nabla u(x)
+\IZ [u(y)-u(x)-\nabla& u(x)\cdot \floo yx ]\, \Jb^P _{ t,x}(dy),\\
\ &(t,x)\in\ioZ, \ \pb\ae
\end{split}
\end{align}
where  $\Jb^P_t$ is the unique solution of  
\begin{align}\label{eq-65p}
\pp_t(dy) \Jb^P _{ t,y}(dx)=\pp_t(x)\Jf^P _{ t,x}(dy),
\end{align}
and the backward drift  $\Bb ^{ P, \delta}$ is given by
\begin{align}\label{eq-67p}
( \Bf ^{ P, \delta}+\Bb ^{ P, \delta})(t,x)
	= \IZ \floo yx \ (\Jf^P _{t,x}+\Jb^P _{ t,x})(dy),
	\qquad (t,x) \ \pb\ae
\end{align}
where the right hand side of this identity is well defined and  $\pb$-integrable.

In particular, if  \eqref{eq-116} holds for every $t_o>0,$ the above results hold for almost every $0<t\le T.$ 
\end{theorem}

\begin{proof}
It is an immediate corollary of Theorem \ref{res-37}, Theorem \ref{res-36} which is  invoked at Lemma \ref{res-39}, and Lemma \ref{res-18}.
\end{proof}

\begin{remark}
With \eqref{eq-114a}, \eqref{eq-67p} gives
\begin{align*}
\Bf ^{ R, \delta}_t(x)+\Bb ^{ P, \delta}_t(x)
	=\IZ \floo yx
\, (\Jf^R _{ t,x}+\Jb ^P _{ t,x})(dy).
\end{align*}
Of course, when $ \delta=0$ we see that
\begin{align*}
\Bb^P=-\Bf^R.
\end{align*}
\end{remark}

Let us introduce the function
 \begin{equation}\label{eq-56}
\hh(a):=\left\{\begin{array}{ll}
a\log a-a+1, & \textrm{if } a>0,\\
1,& \textrm{if }a=0,\\
\infty,& \textrm{if }a<0.
\end{array}\right.
\end{equation}

\begin{corollary}\label{res-44}
Assume that $R\in\MO$ solves
$
R\in \MP _{  \delta}(\Bf ^{ R, \delta},\Jf^R)
$
and satisfies the Hypotheses \ref{ass-08} and \eqref{eq-115}. For any $0\le t_o< T$ and any nonnegative measurable function $j: \ioZZs\to[0, \infty)$ such that
\begin{align}
\sup _{ t_o\le t\le T, x\in\ZZ}\IZ \hh\big(j(t,x,y)\big)\, \Jf^R _{ t,x}(dy)< \infty,
\end{align}
 define 
\begin{align*}
\begin{split}
\Bf ^{ \delta}_t(x)&:=\Bf ^{ R, \delta}_t(x)
	+\IZ \floo yx (j(t,x,y)-1)\, \Jf^R _{ t,x}(dy),\\
\Jf _{ t,x}(dy)&:= j(t,x,y) \ \Jf^R _{ t,x}(dy),
\end{split}
\end{align*}
where $ \delta=0$ if $R$ satisfies \eqref{eq-31a} (in which case $\Bf=\Bf^R$), and $ \delta=1$ otherwise.

Then, the integral in the expression of $\Btf$ is well defined and for any $\pp _{ t_o}\in\PZ$ such that $H(\pp _{ t_o}|\rr _{ t_o})< \infty,$ the martingale problem $\MP _{ t_o, \delta}(\pp _{ t_o},\Btf,\Jf)$ admits a unique solution $P _{ [t_o, T]}$. This means that 
\begin{align*}
P _{ t_o}=\pp _{ t_o},\qquad
\Bf ^{ P, \delta}_t=\Btf_t, \qquad 
\Jf^P_t=\Jf_t,
\qquad t\in\io.
\end{align*}
Furthermore, $H(P _{ [t_o,T]}|R _{ \io})< \infty$, and for almost every $t\in\io$ we have  $\CcZ\subset \dom\Jb_t^P$,  and the time reversal formula \eqref{eq-677}-\eqref{eq-65p}-\eqref{eq-67p} is verified. 
\end{corollary}

\subsection*{Hypotheses. List of properties}
Let us state now the list of properties which enter the set of  Hypotheses \ref{ass-08}.

\par\medskip
\noindent
\textsf{General hypotheses.}	
\begin{itemize}
	\item\emph{Growth of the drift field.} The {drift} field 
	\begin{align}\label{eq-121}
	\Bf ^{\delta}\ \textrm{ is locally bounded.}
	\end{align}

	\item\emph{Integrability of the jump kernel.}  Next estimate  	\begin{align}\label{eq-32} 
	\sup _{ 0\le t\le T, |x|\le \rho }\IZs (|\xi|^2\wedge 1)\, K  _{ t,x}(d\xi)< \infty,
	\quad \forall \rho \ge 0,
	\end{align}
	implies that the sample paths have a finite quadratic variation with finitely many large jumps (with an amplitude larger than 1, say), almost surely.

\item \emph{Local boundedness of the marginals.}\  The Markov  measure $Q\in\MO$  verifies
\begin{align}\label{eq-94}
\qb( \left\{(t,x):|x|\le \rho\right\}) < \infty,\qquad \forall\rho\ge 0.
\end{align}
Of course this holds when $Q$ is a probability measure.
\end{itemize}

\noindent
\textsf{Large jumps.}

\begin{itemize}	
	\item\emph{Large jumps 1.} We define the range of jumps at $x\in\ZZ$ 	by
	\begin{align}\label{eq-76a}
	 \Delta_K(x):=\inf \Big\{ \Delta\ge 0; \sup _{ t\in\ii}K  _{ t,x}(\{\xi\in\Zs; |\xi|\ge \Delta\})=0 \Big\}
	 \in[0, \infty].
	\end{align}
		The hypothesis on $K$ and $Q\in\MO$  is:
	\begin{align}\label{eq-76}
	\begin{split}
	&\Delta_K \textrm{ is a  locally bounded  function and }\\
	&\IiZZs \1 _{ \{ |x|\le \rho, |\xi|\ge 1 \} } \, \qq_t(dx)K _{ t,x}(d\xi)dt< \infty,
	\quad \forall \rho\ge 0.
	\end{split}
	\end{align}	
%\begin{align}\label{eq-109b} \textrm{for some } \delta_o>0, \sup _{(t,x)\in \ii\times B_ \rho}  K _{ t,x} (B_{ \delta_o}^c)< \infty,\qquad \forall \rho\ge 0.\end{align}
	\item \emph{Large jumps 2.}
	Here $Q$ is assumed to be bounded:
 	\begin{align}
	\label{eq-110}
	Q\in\PO \quad \textrm{and}\quad
	\IiZZs  \1 _{ \{|\xi|\ge 1\}}\,  dt\qq_t(dx)K _{ t,x}(d\xi)< \infty.
	\end{align}
\end{itemize}

\noindent
\textsf{Bounded variation sample paths.}	

\begin{itemize}
	\item \emph{Small jumps.}
	\begin{align}\label{eq-31a}
	\sup _{ 0\le t\le T, |x|\le \rho }\IZs \1 _{ \{|\xi|\le 1\}} |\xi|\, K  _{ t,x}(d\xi)< 		\infty,\quad \forall \rho \ge 0.
	\end{align}
	This strengthening  of \eqref{eq-32} implies that the sample paths have a finite variation almost surely.
	
	\item \emph{Continuity of the jump kernel.}
	\begin{align}\label{eq-771}
	(t,x)\mapsto \IZs f(\xi) (1\wedge |\xi|)\,K  _{ t,x}(d\xi)\quad \in C ^{ 0,0} 		_{ t,x}, \quad \forall f\in C_b(\Zs).
	\end{align}
\end{itemize}

\noindent
\textsf{Unbounded variation sample paths.}	
\ The following hypotheses hold under \eqref{eq-32}  but  when  \eqref{eq-31a} fails, that is when
\begin{align}\label{eq-97}
	\sup _{ 0\le t\le T, |x|\le \rho_o }\IZs \1 _{ \{|\xi|\le 1\}} |\xi|\, K  _{ t,x}(d\xi)= \infty,
	\quad \textrm{for some } \rho_o \ge 0.
	\end{align}
Some sample paths (if not all) might have an unbounded variation. In this case we have to specify the  small jumps mechanism.

\begin{itemize}
	\item \emph{Continuity of the jump kernel.}
	\begin{align}\label{eq-772}
	(t,x)\mapsto \IZs f(\xi) (1\wedge |\xi|^2)\,K  _{ t,x}(d\xi)\quad \in C ^{ 0,0} 		_{ t,x}, \quad \forall f\in C_b(\Zs).
	\end{align}
	
	\item\emph{Small jumps}. 
The forward jump kernel $K$ of the process  satisfies
\begin{align}\label{eq-100}
\1 _{ \{|\xi|\le \delta_o(|x|)\}}\  K _{ t,x}(d\xi)
	=\1 _{ \{|\xi|\le \delta_o(|x|)\}}\ k _{ t,x} (\xi)\ [( \phi _{ t,x})\pf \Lambda](d\xi)
\end{align}
for some 
\begin{itemize}
\item
 positive  function  $\delta_o:  [0, \infty)\to (0,1]$, 
\item
  positive function $k:\ii\times\ZZ\times \{\xi\in\Zs: |\xi|\le 1\}\to (0, \infty),$
\item
  measure space $( \mathcal{A}, \Lambda)$ with $ \Lambda$ a positive measure,
 \item
   mapping
$\phi:\iZ\times \mathcal{A} \to\Zs$.
\end{itemize}
We assume that  $k$ is in $C ^{ 0,1,1}(\ii\times\ZZ\times \{\xi\in\Zs: |\xi|\le 1\})$
and satisfies 
\begin{align}\label{eq-119}
\begin{split}
\sup _{t\in\ii, x\in B _{ \rho}} \IZs \1 _{ \{|\xi|\le \delta_o(|x|)\}}
	\big[&k(x,\xi)|\xi|^2 +|\nabla k(x,\xi)|  |\xi|^2\\
	&+ | \nabla_\xi k(x,\xi)| |\xi|^3\big]\ [( \phi _{ t,x})\pf \Lambda](d\xi)
	< \infty,
	\qquad \forall \rho\ge 0. 
\end{split}
\end{align}
We also assume that $\phi$ is measurable and  for all $(t,x, \alpha)\in\iZ\times \mathcal{A}$,
\begin{align}\label{eq-100a}
\phi _{ t, x}( \alpha)=\nabla \psi _{ t, \alpha}(x),
\end{align}
where $ \psi$ is a numerical function  on $\iZ\times \mathcal{A}$ which is $C ^{ 0,2} _{ t,x}$  and satisfies
\begin{align}\label{eq-100b}
&\inf _{t\in\ii,  \alpha\in \mathcal{A}, x\in\ZZ} \nabla^2 \psi _{ t, \alpha}(x)>-\Id,
\\\label{eq-100c}
&\sup _{t\in\ii,  \alpha\in \mathcal{A}, x\in \ZZ}| \nabla \psi _{ t, \alpha}(x)|< \infty,\\
\label{eq-100d}
&\sup _{t\in\ii,  \alpha\in \mathcal{A}, x\in B _{ \rho}} \frac{|\nabla^2 \psi _{ t, \alpha}(x)|}{|\nabla \psi _{ t, \alpha}(x)|}< \infty,
\quad \forall \rho\ge 0.
\end{align}
Here and below,  $\nabla$ stands for the gradient $\nabla_x$ with respect to $x.$
\end{itemize}

\noindent
\textsf{Close to reversibility.}	

	\begin{itemize}
	\item
	Defining $\pi_t(dxdy):=\qq_t(dx) J  _{ t,x}(dy),$ and its symmetrized $			\tilde\pi_t(dxdy):=\pi_t(dydx),$ we assume that  
	\begin{align}\label{eq-81}
	\sup  \left\{
	 \frac{d\tilde \pi_t}{d\pi_t}(x,y); t\in\ii, x: |x|\le \rho, y\in\ZZ
	\right\} < \infty,
	\qquad \forall \rho\ge 0.
	\end{align}
	\end{itemize}

\noindent
\textsf{Control of the marginal flow.}	

	\begin{itemize}

	\item\emph{Positive regular density.}\ 
	There exist $0\le t_o<T$  	 such that $\qq_t(dx)\ll dx$ for all $t\in [t_o,T]$ 		and 
	\begin{align}\label{eq-88}
	(t,x)\mapsto \frac{d\qq_t}{dx}>0
	\end{align}
	is  positive and in $C ^{ 0,1} ([t_o,T]\times\ZZ).$ 
 
		\item 	There exist $0\le t_o<T$ such that

	\begin{align}\label{eq-89}
	\sup _{t\in[t_o,T], x: |x|\le \rho, \xi\in\supp K _{ t,x}}
	\frac{d \qq^\xi_t}{d\qq_t}(x)< \infty,
	\qquad \forall \rho\ge 0,
	\end{align}
	 where $\qq^\xi:=\tau^\xi\pf\qq$ with  $\tau^\xi(x):=x+\xi,$ $x\in\ZZ,$  the translation by the jump $\xi\in\Zs$.
	
	\end{itemize}

\begin{remarks}[About the Hypotheses \ref{ass-08}]\   \label{rem-01}
\begin{enumerate}[(i)]

\item\emph{About \eqref{eq-100c}.}\  This requirement is not a restriction because $\phi=\nabla \psi$ represents a \emph{small} jump. 

\item\emph{About \eqref{eq-100a} -- \eqref{eq-100d}.}\  One of the simplest functions $ \psi$ satisfying \eqref{eq-100b}, \eqref{eq-100c} and \eqref{eq-100d} is 
\begin{align*}
\psi _{ t, \alpha}(x)=  \alpha\scal x
\end{align*}
with $ \alpha\in \mathcal{A}= \left\{ \alpha\in\Zs, | \alpha|\le 1\right\} .$ Choosing $\Lambda( d \alpha)=d \alpha$  gives $[( \phi _{ t,x})\pf \Lambda](d\xi)=d\xi$ and \eqref{eq-100} is
$
K _{ t,x}(d\xi)= k _{ t,x} (\xi)\,d\xi
$
in restriction to $|\xi|\le \delta_o(|x|).$

\item\emph{About \eqref{eq-81} and \eqref{eq-89}.}\  

\begin{enumerate}[-]
\item
Under the Hypotheses \ref{ass-08}-a-b-c, one of the assumption \ref{ass-08}-d1 or \ref{ass-08}-d2, i.e.\   \eqref{eq-81} or \eqref{eq-89}, implies \ref{ass-08}-d3, i.e.\ $\LLb u$ is bounded, see Lemmas \ref{res-39} and \ref{res-42}.

\item\emph{About \eqref{eq-81}.}\  
If the   path measure is reversible  then   $\pi=\tilde \pi,$ in which case the quantity in \eqref{eq-81} is equal to $1$ for all $ \rho$. Hence  this quantity measures some proximity to reversibility  of the  jump mechanism.

\item\emph{About \eqref{eq-89}.}\  
It is proved at Proposition  \ref{res-33} that the time reversal formula implies that
$\qq^\xi\ll\qq$ for any "admissible" jump $\xi.$ The only restriction in this hypothesis is the local boundedness of the derivative $d\qq^\xi/d\qq.$
\end{enumerate}
\end{enumerate}
\end{remarks}

\section{Examples}
\label{sec-examples}

We look at families of examples where no assumption is required on the marginal flow.

\subsection*{Example 1}

The easiest setting corresponds to $P\in\PO$ with bounded variation sample paths.  The assumptions are
\begin{enumerate}[-]
\item
the drift field $\Bf\in C ^{ 0,1} _{ t,x}$ is such that there exists some $c\ge0$ such that 
$\Bf_t(x)\cdot  x\le c (1+|x|^2)$ for all $(t,x)\in\iZ.$
\item
the jump kernel $\Kf$ allows for bounded variation sample paths:
\begin{align*}
	\sup _{ 0\le t\le T, |x|\le \rho }\IZs (|\xi|\wedge 1)\, \Kf  _{ t,x}(d\xi)< \infty,
	\quad \forall \rho \ge 0,
	\end{align*}   
	its range $ \Delta _{ \Kf}$ is locally bounded, recall \eqref{eq-76a}, and it is continuous in the sense of \eqref{eq-771}.
\end{enumerate}

\begin{proposition}
Under these assumptions, for any initial distribution $\pp_0\in\PZ,$ there exists a unique solution  $P\in\MP(\pp_0,\Bf,\Kf)$. 
\\
For almost every $t\in\ii$ we have  $\CcZ\subset \dom\LLb_t^P$ and   
 the time reversal formula \eqref{eq-677}-\eqref{eq-65p}-\eqref{eq-67p} is verified.  
\end{proposition}

\begin{proof}
The existence result is standard, and the time reversal statement is a direct corollary of Theorem \ref{res-36}.
\end{proof}

Extending this result with the entropic improvement (Theorem \ref{res-43}) requires a control of the marginal flow as in  \eqref{eq-89}. As shown by next illustration, this control (if available) must be done on $[t_o,T]$ for any $t_o>0$. But this is enough to recover the time reversal formula on $(0,T].$

\subsubsection*{Poisson process with parameter $\lambda>0$} This corresponds to the state space $ \mathbb{N},$ $\pp_0= \delta_0,$ $\Bf=0$ and $\Kf _{ t,x}= \lambda \delta_1$ for all $t,x.$
 As $\pp_t(k)= e ^{ - \lambda t} ( \lambda t)^k/k!,$ we see that 
 the time reversal formula $
\pp _{t}(k)\Jb _{ t,k}(k-1)=\pp_t(k-1)\Jf _{ k-1}(k)
$	
implies that
\begin{align*}
\Jb _{ t,k}(k-1)= \frac{ \lambda \pp_t(k-1)}{\pp _{t}(k)}
	=k/t,
	\qquad k\ge 1.
\end{align*}
Remark that it does not depend on $ \lambda.$
On the other hand, the density in \eqref{eq-89}:
\begin{align*}
\frac{d \tau ^{ 1}\pf\pp_t}{d\pp_t}(k)= \frac{e ^{ - \lambda t} ( \lambda t) ^{ k-1}/(k-1)!}{e ^{ - \lambda t} ( \lambda t) ^{ k}/k!}
	= \frac{k}{ \lambda t} ,
\end{align*}
explodes as $t$ tends to zero.

\subsection*{Example 2}
We prove the time reversal formula for some Markov measure without drift whose jump kernel   is absolutely continuous. More precisely, 
\begin{align*}
\Jf_{ t,x}(dy)
	= e ^{ -[V(y)-V(x)]/2}\ j(t,x,y) \sigma(|y-x|)\ dy,
\end{align*}
where it is assumed that 
\begin{align*}
V:\ZZ\to\RR \textrm{ is }C^1 \textrm{ and } \inf V> - \infty,
\end{align*}
and 
 $ \sigma: (0, \infty)\to[0, \infty)$ is  a continuous  function which is differentiable  on $(0,r_o)$ for some $r_o>0$  satisfying 
\begin{align}\label{eq-124}
\begin{split}
&(a)\quad
\sup _{ x:|x|\le \rho}\IZs  \1 _{ \{|\xi|\ge r_o\}} e ^{ - V(x+\xi)/2} \sigma(|\xi|)\ d\xi< \infty,
\qquad \forall \rho\ge 0,
\\
&(b)\quad 
\IZZs  \1 _{ \{|\xi|\ge r_o\}} e ^{ -[V(x)+V(x+\xi)]/2} \sigma(|\xi|)\ dxd\xi< \infty,\\
&(c)\quad 
\int _{ (0,r_o)} \sigma(r) r ^{ n+1}\, dr< \infty,\\ 
&(d)\quad 
\int _{ (0,r_o)} \sigma'(r) r ^{ n+2}\, dr< \infty,
\end{split}
\end{align}
and the initial marginal  $\pp_0\in\PZ$ is absolutely continuous and satisfies
\begin{align}\label{eq-126}
\IZ \1 _{ \{d\pp_0/dx >0\}} \Big(\log \big( \frac{d\pp_0}{dx}) +V(x)\Big)\, \pp_0(dx)< \infty.
\end{align}
The measurable function
 $j: \iZZ\to[0, \infty)$ is   such that
\begin{align*}
\sup _{0\le t\le T, x\in\ZZ}\IZ \hh\big(j(t,x,y)\big)\, e ^{ -[V(y)-V(x)]/2}\ \sigma(|y-x|)\ dy< \infty,
\end{align*}
or more generally
\begin{align}\label{eq-127}
\IZ \hh\big(j(t,x,y)\big)\, e ^{ -[V(y)-V(x)]/2}\ \sigma(|y-x|)\ dt\pp_t(dx)dy< \infty,
\end{align}
where $\hh$ is defined at \eqref{eq-56}.

\begin{proposition}\label{res-47}
The martingale problem $\MP _{ \delta}(\pp_0,0,\Jf)$ admits a unique solution $P\in\PO$. It satisfies $H(P|R)< \infty,$ and more precisely $H(P|R) < \infty$ if and only if \eqref{eq-126} and \eqref{eq-127} are satisfied. 
\\
For almost every $t\in\ii$ we have  $\CcZ\subset \dom\LLb_t^P$ and   
 the time reversal formula \eqref{eq-677}-\eqref{eq-65p}-\eqref{eq-67p} is verified.  
\end{proposition}

\begin{proof}
It is a direct consequence of Corollary \ref{res-44} with Lemma \ref{res-46} below. There exists  a unique solution to the martingale problem because our assumptions imply that $H(P|R)< \infty$ and a fortiori that $P\ll R,$ where $R$ is the unique solution of its martingale problem by Lemma  \ref{res-46}. See \cite{Jac79} for this powerful argument.
\end{proof}

\subsubsection*{A reversible jump process}

Consider  the equilibrium measure $\rr\in\MZ$ defined by
\begin{align*}
\rr(dx)
	:= e ^{ -V(x)}\ dx,
\end{align*}
and the jump kernel 
\begin{align*}
J_x^R(dy)
	:=  e ^{ -[V(y)-V(x)]/2}\ s(x,y)\ dy
\end{align*}
where $s$ is a nonnegative measurable \emph{symmetric} ($s(x,y)=s(y,x)$) function defined on $\ZZZ$ minus its diagonal,  such that 
\begin{align}\label{eq-118}
\int _{ \{x\neq y\}} (1\wedge |y-x|^2)\ e ^{ -[V(y)-V(x)]/2}\ s(x,y)\ dxdy< \infty,
\end{align}
and for any $ \rho\ge 0,$ 
\begin{align}\label{eq-117}
\begin{split}
&\sup _{ x:|x|\le \rho}\IZ \1_ {\{ 0<|y-x|\le 1\}} |y-x|^2\ s(x,y)\ dy< \infty,\\
&\sup _{ x:|x|\le \rho} \frac{s(x-\xi,x)}{s(x,x+\xi)}=1+O _{ \xi\to 0}(|\xi|),
\qquad \xi\in\Zs.
\end{split}
\end{align}
For any $u\in\CcZ$ and any $x\in\ZZ,$
\begin{align*}
\LLf^Ru(x)
	=\IZ [u(x+\xi)-u(x) ]  e ^{ -[V(x+\xi)-V(x)]/2}\ s(x,x+\xi)\ d\xi
\end{align*}
is well-defined with an abuse of notation but without introducing any truncation $\xd$. To see that this is true, control the small jump contribution by writing the integral with respect to $\xi$ as its half sum  with the same integral  after the change of variables $\xi\to -\xi,$ use the symmetry of $s$ and conclude with \eqref{eq-117}.

\begin{proposition}\label{res-45}
If the solution of the martingale problem 
$
R\in\MP(\rr,0,J^R)
$
 exists, then it  is reversible,  that is: $$\LLb^R=\LLf^R.$$ 
\end{proposition}

\begin{proof}
A direct computation using the symmetry of $s$   shows that the formal adjoint (roughly speaking: in $L^2(\Leb)$) $\LLf ^{ R, \ast}$ of $\LLf^R$  annihilates $e ^{ -V}:$ $\LLf ^{ R, \ast} e ^{ -V}=0.$ This proves that $\rr$ is a stationary measure.
\\
On the other hand, the estimate \eqref{eq-118} implies the hypothesis \eqref{eq-84}-b 
 of Theorem \ref{res-37} which tells us that it remains to verify that the flux equation \eqref{eq-FE} is valid in this situation, i.e.\ 
$	%\begin{align*}
\rr(dy) J_y^R(dx)=\rr(dx) J_x^R(dy).
$	%\end{align*}
But this amounts to
\begin{equation*}
e ^{ -[V(y)+V(x)]/2}\ s(y,x)\ dxdy
	= e ^{ -[V(y)+V(x)]/2}\  s(x,y)\ dxdy,
\end{equation*}
and is true because $s$ is symmetric.
\end{proof}

Let us go back to a function  $s$ of the form $s(x,y)=  \sigma(|y-x|).$

\begin{lemma}\label{res-46}
The $\rr$-reversible path measure $R\in\MP(\rr,0,J^R)$  exists and it satisfies   \eqref{eq-115} and the Hypotheses \ref{ass-08}-a-b2-c-d1.
\end{lemma}

\begin{proof}
The jump kernel writes as
\begin{align*}
K^R_x(d\xi)= k(x,\xi)\, d\xi,
\qquad k(x,\xi):=e ^{ -[V(x+\xi)-V(x)]/2} \sigma(|\xi|).
\end{align*}
The existence of a solution to the martingale problem follows from the existence of a unique strong solution due to the regularity and integrability of the the kernel, plus Yamada's theorem, see \cite{Jac79}. Furthermore  in this case \eqref{eq-115} holds trivially.

The first  requirement of \eqref{eq-117} becomes
$
\IZs \1_ {\{|\xi|\le 1\}} |\xi|^2\ \sigma(|\xi|)\, d\xi< \infty,
$
which amounts to \eqref{eq-124}-(c) 
and the second one is trivially satisfied. The control of large jumps is done by \eqref{eq-110} which is finite by \eqref{eq-124}-(b).
\\
 In view of Remark \ref{rem-01}-(iii), to control the small jumps it remains to verify \eqref{eq-119}:
\begin{align*}
\sup _{x\in B _{ \rho}} \IZs \1 _{ \{|\xi|\le 1\}}
	\big[&k(x,\xi)|\xi|^2 +|\nabla k(x,\xi)|  |\xi|^2
	+ | \nabla_\xi k(x,\xi)| |\xi|^3\big]\, d\xi
	< \infty,
	\qquad \forall \rho\ge 0. 
\end{align*}
Because $V$ is $C^1,$ the contribution of the first two terms in the integrand is finite by \eqref{eq-124}-(c).   Similarly, the contribution of the third term is finite by 
\eqref{eq-124}-(d). Finally, Hypothesis \ref{ass-08}-d1 is trivially satisfied because $R$ is reversible by Proposition \ref{res-45}.
 \end{proof}

The proof of Proposition \ref{res-47} is twofold: (1) obtain $\LLb^R=\LLf^R$ for a reversible reference Markov measure $R,$ (2) then extend the time reversal formula to $P$ such that $H(P|R)< \infty.$  The main advantage of this strategy is that \emph{we do not have to suppose any a priori regularity of the marginals of $P,$ } and this works for a vast family of Markov measures.

\subsection*{Example 3}
A limitation of  Example 2 is that the support of the jump kernel $K^R$ is symmetric, i.e.\ $\xi\in\supp K^R\Leftrightarrow -\xi\in\supp K^R.$  Of course, the density $j$ may vanish at some places, allowing for asymmetric jumps for $P$. But the entropic price to pay for such a killing is 
$	%\begin{align*}
\IiZZ \1 _{ \{j(t,x,y)=0\}} e ^{ -[V(y)-V(x)]/2} \sigma(|y-x|) \pp_t(dx)dydt,
$	%\end{align*}
which might be infinite if there are too many small jumps.  
  In this subsection, we look at  examples with not necessarily diffuse   kernels  and possibly asymmetric small jumps. The jump kernel is 
\begin{align*}
\Kf _{ t,x}(d\xi)
	=k(t,x,\xi)\ \Lambda(d\xi)
\end{align*}
where the nonnegative measure  $ \Lambda$  on $\Zs$ verifies
\begin{align*}
\IZs (1\wedge|\xi|^2)\, \Lambda(d\xi)< \infty.
\end{align*}
and the measurable function
 $k: \iZZs\to[0, \infty)$ is   such that
\begin{align*}
\sup _{0\le t\le T, x\in\ZZ}\IZs \hh\big(k(t,x,\xi)\big)\, \Lambda(d\xi)< \infty,
\end{align*}
or more generally
\begin{align}\label{eq-128}
\IiZZs \hh\big(k(t,x,\xi)\big)\, dt\pp_t(dx)\Lambda(d\xi)< \infty,
\end{align}
where $\hh$ is defined at \eqref{eq-56}.
\\
The initial marginal  $\pp_0\in\PZ$ is absolutely continuous and satisfies
\begin{align}\label{eq-129}
\IZ \1 _{ \{d\pp_0/dx >0\}} \log \big( \frac{d\pp_0}{dx} \big)\, \pp_0(dx)< \infty.
\end{align}

\begin{proposition}\label{res-exple3}
The martingale problem $\MP _{ \delta}(\pp_0,0,\Kf)$ admits a unique solution $P\in\PO$. We have  $H(P|R)< \infty,$ and more precisely $H(P|R) < \infty$ if and only if \eqref{eq-128} and \eqref{eq-129} are satisfied. 
\\
Furthermore,  $\CcZ\subset \dom\LLb_t^P$ for almost every $t\in\ii$,    and   
 the time reversal formula \eqref{eq-677}-\eqref{eq-65p}-\eqref{eq-67p} is verified.  

\end{proposition}

\begin{proof}
It is a  consequence of Theorem \ref{res-43} with Lemma \ref{res-48} below.
This lemma states that  $R\in\MO$ is the law of a stationary process with independent increments, and the unique solution to its martingale problem.    uniqueness and \eqref{eq-115}. Note that although $R$ do not meet the hypothesis \ref{ass-08}-b about the large jumps, the result  of Theorem  \ref{res-43} is still valid because a direct inspection shows that $\LL^Ru$ is bounded. Consequently we do not need Lemmas \ref{res-39} and \ref{res-42} which require this assumption and whose purpose is to obtain this boundedness.
\\
The path measure $P$ is also  the unique solution to its martingale problem because our assumptions imply that $H(P|R)< \infty$ and a fortiori that $P\ll R,$ where $R$ is the unique solution of its martingale problem. 
\end{proof}

\subsubsection*{A stationary jump process}

The  forward generator of the reference measure $R$ is defined for any $u\in\CcZ$ and $x\in\ZZ$ by
\begin{align*}
\LLf^Ru(x)=b\scal \nabla u(x)+\IZs [u(x+\xi)-u(x)-\nabla u(x)\cdot \flo \xi]\, \Lambda(d\xi)
\end{align*}
with $b\in\ZZ$ a vector and a jump kernel 
\begin{align*}
K^R _{ t,x}(d\xi)= \Lambda(d\xi),
\quad \forall (t,x)\in\iZ,
\end{align*}
which do not depend on $(t,x).$ This is the generator of a process with stationary independent increments.  The initial marginal of $R$ is Lebesgue measure: $R_0=\Leb$.

\begin{lemma}\label{res-48}
The path measure $R\in\MO$ is the unique solution of $\MP_ \delta(\Leb,b, \Lambda).$ It is the law of a $\Leb$-stationary process. Its time reversal $R^*\in\MO$ is also the law of a $\Leb$-stationary process with independent stationary increments, and is  the unique solution of $\MP_ \delta(\Leb,-b, \Lambda^*)$ with
\begin{align*}
\Lambda^*= \eta\pf \Lambda,
\end{align*}
where $\eta(\xi):=-\xi,$ $\xi\in\Zs.$
\end{lemma}

\begin{proof}
As in Proposition \ref{res-45}'s proof, we rely on Theorem \ref{res-37}. 
\\
Suppose for a while that $ \Lambda$ has a bounded support.  Using this  assumption, in particular to show that for any  $u\in\CcZ$, $\LLf^R u\in L^1\cap L^2(\Leb),$ a direct computation  shows that  the adjoint of $\LLf^R$ is given by
 \begin{align*}
\LLf ^{ R,\ast}u(x)=-b\scal \nabla u(x)+\IZs [u(x+\xi)-u(x)-\nabla u(x)\cdot \flo \xi]\, \Lambda^*(d\xi).
\end{align*}
By a standard argument, this shows that $R$ is $\Leb$-stationary. We are now in position to apply  Theorem \ref{res-37} under the assumption \eqref{eq-84}-a,  with $\qq=\Leb$.
Therefore,  $R^*\in\MP _{  \delta}(\Leb,-b, \Lambda^*).$

To extend the result to the case where $ \Lambda$ has an unbounded support, consider for any $k\ge 1$ its restriction $ \Lambda^k:= \1 _{ B_k}\, \Lambda$ to the ball  of radius $k$ in $\Zs$. We have just shown that  $R^k\in\MP_\delta(\Leb,b, \Lambda^k)$ satisfies $R ^{ k,\ast}\in\MP_\delta(\Leb,-b, \Lambda ^{ k,\ast})$ with $\Lambda ^{ k,\ast}:=\eta\pf \Lambda^k=\1 _{ B_k} \Lambda^*.$  Clearly $\Lim k \Lambda ^{ k,\ast}= \Lambda^*$ for the weak topology defined by the continuous test functions $ \varphi$ on $\Zs$ such that $\sup _{ \xi\in\Zs} | \varphi(\xi)|/(1\wedge|\xi|^2)< \infty.$ It follows that $\Lim k R ^{ k,\ast}\in\MP_\delta(\Leb, -b, \Lambda^*)$ for the narrow topology, see \cite{JaShi87}.\\
On the other hand, $\Lim k R ^{ k,\ast}=R^*$, because  for the same reason $\Lim k R^k=R,$ and time reversal is continuous. We have proved that 
  $R^*\in\MP_ \delta(-b,0, \Lambda^*)$, as announced.
\end{proof}

Another expression of Lemma \ref{res-48} is
\begin{align*}
\LLb ^{ R}u(x)=-b\scal \nabla u(x)+\IZs [u(x+\xi)-u(x)-\nabla u(x)\cdot \flo \xi]\, \Lambda^*(d\xi).
\end{align*}

\subsection*{Adding a drift term}

This  is a remark about  a remaining difficulty which is not overcome in this paper.  Adding a drift term, that is considering  $\MP _{ \delta}(\pp_0,b,\Kf)$ instead of  $\MP _{ \delta}(\pp_0,0,\Kf)$ in Proposition \ref{res-exple3},  would require  to adapt standard proofs of existence of flows of diffeomorphisms  so that one can incorporate random jumps, and  also to  build  synchronous couplings to try to obtain some control of the regularity of the time marginals. In any case, this does not seem to be an easy improvement to achieve.

\section{Integration by parts formula}
\label{sec-IbP}

Our time reversal results  rely on an integration by parts (IbP) formula  for the carré du champ which was proved by  Cattiaux, Gentil and the authors in \cite{CCGL20}.
Before stating this IbP formula  at Theorem \ref{res-02}, let us recall the definitions of Markov measures and extended generators.

\subsection*{Carré du champ}

Let $Q$ be a path measure on $\OO.$
Its forward  carré du champ  is the forward-adapted process defined by
\begin{align*}
\Gaf^Q_t(u,v):= \LLf^Q_t(uv)-u\LLf^Q _tv-v\LLf^Q_t u,
\quad (u,v)\in\dom \Gaf^Q_t,\ 0\le t\le T,
\end{align*}
where 
$
\dom\Gaf^Q_t:= \left\{(u,v); \  u, v, uv \in\dom \LLf^Q_t \right\} .
$
\\
We introduce a class $\UU$ of functions on $\XX$ such that
\begin{align*}
 \mathcal{U}\subset \dom \LLf^Q_t\cap C_b(\XX)
\end{align*}
for all $0\le t\le T$ and any path measure $Q$ of interest, where $C_b(\XX)$ is the space of all bounded continuous functions on $\XX.$  \emph{We assume that $\UU$ is an algebra,} i.e.\ $u,v\in\UU$ implies $uv\in\UU.$ 
In particular, 
\begin{align}\label{eq-43}
u,v\in\UU\implies (u,v)\in\dom \Gaf^Q _t.
\end{align}
We shall mainly consider functions in $\UU$ and make an intensive use of their carré du champ. In each setting, this algebra will be chosen rich enough to determine a Markov dynamics, i.e.\ to solve in a unique way some relevant  martingale problem. We shall see that $\UU=\CcZ $ is a good choice.

\begin{remark}\label{rem-04}
The requirement that {$ \mathcal{U}$ is an algebra} (it is necessary that $uv$ belongs to $\dom \LLf^Q$ to consider $\LLf^Q (uv)$), is  strong.  Let us say that a semimartingale whose bounded variation term is absolutely continuous is ``nice''. The product of two semimartingales is a semimartingale, but  the product of two { nice} semimartingales might not be  nice anymore.   In general, a martingale representation theorem is needed to verify the stability of the product of nice semimartingales.
\end{remark}

Next result is the cornerstone of the proofs of  time reversal formulas.
We introduce the class of functions
\begin{equation}\label{eq-42}
\UQ := \Big\{ u\in \UU; \LLf^Qu\in L^1(\qb),\ \Gaf^Q (u)\in L^1(\qb)  \Big\}.
\end{equation}

\begin{theorem}[IbP of the carré du champ, \cite{CCGL20}]\label{res-02}
Let $Q\in\MO$ be any Markov measure. Take two functions $u,v$ in  $\UQ.$

\begin{enumerate}[(a)]
\item
If
\begin{align}\label{eq-39a}
u\in\dom \LLb^Q\quad \textrm{and}\quad \LLb^Qu\in L^1(\qb),
 \end{align} 
then for almost every $t$ 
\begin{align}\label{eq-24}
\IZ \Big\{(\LLf_t^Qu+\LLb_t^Qu)\, v + \Gaf_t^Q(u,v) \Big\}\, d\qq_t
=0.
\end{align}

\item
Suppose that  
\begin{equation}\label{eq-39b}
(t,x)\mapsto \Gaf^Q_t(u,v)(x) \textrm{ is continuous},
\end{equation}
the class of functions $\UQ$ determines the  weak convergence of Borel measures on $\XX$,
and the linear form 
\begin{equation}\label{eq-39c}
w\in \UU ^{ \Qb} \mapsto \IiZ \Gaf^Q_t(u,w_t)(x) \,dt \qq_t(dx)
\end{equation}
on $\UU ^{ \Qb}:= \left\{ w\in C_b(\iX); \ w(t,\sbt)\in \UQ,\ \forall 0\le t\le T\right\} $ defines a finite measure on $\iX$.\\
Then, \eqref{eq-39a} and \eqref{eq-24} are satisfied.
\end{enumerate}
\end{theorem}

\begin{remarks}\ 
\begin{enumerate}[(a)]
\item
Statement  \eqref{eq-39c} is an IbP formula for the carré du champ when $\Gaf$ is regular, while \eqref{eq-24} extends it to a non-regular setting under the weaker  condition \eqref{eq-39a}.

\item
The backward carré du champ  is defined by
\begin{align*}
\Gab^Q_t(u,v):= \LLb^Q_t(uv)-u\LLb^Q_t v-v\LLb^Q_t u,
\end{align*}
for any $0\le t\le T$ and $(u,v)\in\dom \Gab^Q_t.$ As we shall see, in general for a jump path measure $Q$, contrary to continuous diffusion processes, we have:  $\Gaf^Q\neq\Gab^Q.$ However,  we see with the IbP formula \eqref{eq-24} that: $\IZ \Gaf^Q(u,v)\,d\qq_t=\IZ \Gab^Q(u,v)\,d\qq_t,$ for $u,v$  in $\UQ,$ as soon as $u$ verifies \eqref{eq-39a}.
\end{enumerate}
\end{remarks}

\section{Abstract characterization} \label{sec-jp-abs}

In this section, the IbP formula of Theorem \ref{res-02} is used to obtain at Theorem \ref{res-37} an abstract characterization in a   general setting for the validity of  a time reversal formula for  a Markov jump process on $\XX=\ZZ$.   In next Sections \ref{sec-regular-kernel} and \ref{sec-entropy}, we work out explicit  assumptions  which verify this criterion, and  therefore warrant the time reversal formula.

\subsection*{Jump process on $\ZZ$}
Let us recall basic notions about jump processes on $\ZZ.$

\subsubsection*{Test functions}
By Itô's formula, under our boundedness hypotheses, for any Markov measure $Q\in\MO$ with   generator \eqref{eq-60}, we have: $\CcZ  \subset \dom\LL^Q$ and $\LL^Q=\LL$ in restriction to $\CcZ  .$ This is a good reason for choosing 
$
\mathcal{U}=\CcZ .
$

 \subsubsection*{Carré du champ}

The  carré du champ is 
 \begin{align*}
 \Gamma_t(u,v)(x)= \IZs [u(x+\xi)-u(x)][v(x+\xi)-v(x)]\, K _{ t,x}(d\xi),
 \quad u,v\in \CcZ  .
 \end{align*}
Remark that $ \mathcal{U}=\CcZ $ is an algebra, as required by the hypotheses of  the IbP formula.

\subsection*{Statement of the time reversal formula}

The time reversal formula is easier to grasp when written with $J$ rather than $K.$
The path measure $Q\in\MO$ is such that for any
 function $u$ in $\UU$, 
\begin{align}\label{eq-752}
\LLf^Q_t u(x)= \Btf_t(x)\cdot \nabla u(x)+\IZ [u(y)-u(x)-\nabla u(x)\cdot \floo yx ]\, \Jf^Q_{ t,x}(dy).
\end{align}
In the simple case where \eqref{eq-871} holds, the forward generator $\LLf^Q$ is \eqref{eq-60b}, that is
\begin{align}\label{eq-751}
\LLf^Q_t u(x)= \Bf_t(x)\cdot \nabla u(x)+\IZ [u(y)-u(x)]\, \Jf^Q _{ t,x}(dy).
\end{align}

\subsubsection*{Flux equation}

In analogy with \eqref{eq-65p}, we introduce the flux equation $\FE(\qq,\Jf):$
\begin{align}\label{eq-FE}
\qq(dy) \mathcal{J} _{ y}(dx)=\qq(x)\Jf_x(dy),
\end{align}
where the unknown is the kernel $ \mathcal{J}=\{ \mathcal{J}_y, y\in\ZZ\}$  and the known coefficients are the  positive measure $\qq$ and the forward kernel $\Jf.$  
We shall see at Theorem \ref{res-37} that under some additional hypotheses, $Q$ admits a time reversal formula if and only if the equation $\FE(\qq_t,\Jf^Q_t)$ admits a solution for almost every $t$. By next Proposition \ref{res-33}, if it exists, this solution is unique. The jump kernel part of the time reversal formula states precisely that the backward kernel $\Jb^Q_t$ is the solution of  $\FE(\qq_t,\Jf^Q_t)$.

 Next proposition   gives an if-and-only-if  condition for the existence of a solution to $\FE(\qq,\Jf)$ and asserts its (already announced)  uniqueness.
\\
For any measurable nonnegative  function  $ \sigma:\ZZZ\mapsto [0, \infty)$  and any $\qq\in\MZ,$ define the measure $\If{\qq}{ \sigma}$ on $\ZZ$  by
\begin{align*}
\If{\qq}{ \sigma}(\sbt):=\IZ [\sigma\Jf] _{ x}(\sbt)\,\qq(dx),
\end{align*}
where $ [\sigma\Jf]_{ x}(dy):= \sigma(x,y)\Jf _{ x}(dy).$

\begin{proposition}\label{res-33}
The equation \eqref{eq-FE}: $\FE(\qq,\Jf)$ admits a solution  if and only if 
\begin{align} 
\label{eq-70}
\If{\qq}{ \sigma}\ll\qq,
\end{align}
for some  measurable \emph{positive}  function  $ \sigma:\ZZZ\mapsto (0, \infty)$ such that
$\If{\qq} \sigma(\ZZ) < \infty.$
\\
If this holds for one function $ \sigma$, then it holds for all measurable positive function $ \sigma$ satisfying  $\If{\qq} \sigma(\ZZ) < \infty.$
\\
Moreover,  the solution $\mathcal{J}$ of $\FE(\qq,\Jf)$ is unique and 
\begin{equation} \label{eq-83}
\begin{split}
\IZ \sigma(x,y)\, \mathcal{J} _{ y}(dx)= \frac{d\If{\qq} \sigma}{d\qq}&(y) \\
	&= \IZs  \frac{d \qq^\xi}{d\qq}(y)\, \sigma(y-\xi, y)\, \Kf_{ y-\xi}(d\xi),
\qquad \forall y\ \qq\ae
\end{split}
\end{equation}
where $\qq^\xi:=\tau^\xi\pf\qq$ with  $\tau^\xi(x):=x+\xi,$ $x\in\ZZ,$  the translation by the jump $\xi\in\Zs$. The identity  \eqref{eq-83} is valid even if $ \sigma$ vanishes at some places.
\end{proposition}

It is part of the result that, when \eqref{eq-70} is satisfied, 
 the Radon-Nikodym derivative  $[d \qq^\xi/d\qq](y)$ is well-defined for almost every $(y,\xi)$ with respect to the measure
$\qq \sigma\Kf(dy,d\xi):= \qq(dy) \sigma(y-\xi,y)\Kf _{ y-\xi}(d\xi).$

Next theorem is the main result of this section.

\begin{theorem}[Time reversal formula]\label{res-37}
Under the Hypothesis \ref{ass-08}-a, suppose that  \eqref{eq-872} is replaced by the stronger requirement 
\begin{align}\label{eq-84}
\begin{split}
&(a)\quad \eqref{eq-76} \quad \textrm{and}\quad \IiZZs \1 _{ \{ |x|\le \rho \} } (1\wedge |\xi|^2)\, \qq_t(dx)\Kf^Q _{ t,x}(d\xi)dt< \infty,
\quad \forall \rho\ge 0,\\&\textrm{or}\\
&(b)\quad \IiZZs   (1\wedge |\xi|^2)\, \qq_t(dx)\Kf^Q _{ t,x}(d\xi)dt< \infty.
\end{split}
\end{align}
Assume also that $Q\in\MO$ is such that
\begin{align*}
C_c^2(\ZZ)=:\UU=\UQ.
\end{align*}
\begin{enumerate}[(a)]
\item
Then, any $u\in C_c^2(\ZZ)\subset \dom \LLf^Q$ is such that $\LLf^Qu$ is $\qb$-integrable.

 \item
Moreover, the three following statements are equivalent:
\begin{align}\label{eq-86}
(1)\qquad &u\in\dom\LLb^Q,\quad \LLb^Qu\in L^1(\qb),\qquad \forall u\in C_c^2(\ZZ);\\
\label{eq-93}
(2)\qquad & \FE(\qq_t,\Jf^Q_t), \textrm{ see  \eqref{eq-FE}, admits a solution for almost every }t\in\ii; \\
\label{eq-70b}
(3)\qquad &\If{t,\qq_t}{ \sigma}\ll\qq_t,
\quad \textrm{for almost every }t,\\ 
\nonumber 
&\textrm{where $ \sigma>0$ verifies }\IiZZ \sigma(x,y)\, \qq_t(dx)\Jf^Q _{ t,x}(dy)dt< \infty.
\end{align}
The estimate \eqref{eq-84} implies that  the function $ \sigma$ in \eqref{eq-70b} can be chosen of the form  $ \sigma(x,y)=\tilde  \sigma(x)(1\wedge | y-x|^2)$,   for some positive function  $\tilde \sigma:\ZZ\to(0, \infty)$.

\item
In this case, $\LLb^Q$ is given by 
\begin{align}\label{eq-85}
\LLb^Q_t u(x)= \Btb_t(x)\cdot \nabla u(x)+\IZ [u(y)-u(x)-\nabla u(x)\cdot \floo yx ]\, \Jb^Q _{ t,x}(dy),
\end{align}
where  for almost every $t$, $\Jb^Q_t$ is the unique solution of $\FE(\qq_t,\Jf^Q_t)$, that is
\begin{align}\label{eq-65}
\qq_t(dy) \Jb _{ t,y}(dx)=\qq_t(x)\Jf _{ t,x}(dy),
\end{align}
and the backward drift  $\Btb$ is given by 
\begin{align}\label{eq-67}
( \Btf+\Btb)(t,x)
	= \IZ \floo yx \ (\Jf^Q _{t,x}+\Jb^Q _{ t,x})(dy),
	\quad \qb\ae
\end{align}
where the right hand side of this identity is well defined and  $\qb$-integrable.
\end{enumerate}
 \end{theorem}

\begin{remarks}\ 
\begin{enumerate}[(a)]
\item
Roughly speaking, \eqref{eq-70b} implies that 
\begin{align*}
\supp(\Jf^Q _{ t,x})\subset \supp(\qq_t),
\quad \forall (t,x)\ \qb\ae
\end{align*}
and also that if $\qq_t$ is a diffuse measure on some subset $A\subset \supp(\qq_t),$ the jump mechanism is not allowed to ``create'' singular structures such as ``Dirac or Cantor masses''   in $A.$ It is likely that this must hold for a large class of non-pathological Markov processes. 

\item
When the sample paths have bounded variations, i.e.\ under hypothesis \eqref{eq-871}, one chooses $\floo yx =0,$ and writes $\Bf$ and $ \Bb$ without delta, so that \eqref{eq-67} is simply  $$\Bf+\Bb=0.$$
In the general case, this still holds in average : 
\begin{align*}
\IZ (\Btf_t+\Btb_t)\, d\qq_t=0
\end{align*}
for almost all                         
 $t$, because 
 \begin{multline*}
 \IZ (\Btf_t+\Btb_t)\, d\qq_t
 	= \IZ \floo yx  \,\qq_t(dx) [\Jf _{ t,x}(dy)+\Jb _{ t,x}(dy)]\\
	= \IZ \floo yx  \,[\qq_t(dx) \Jf _{ t,x}(dy)+\qq_t(dy)\Jf _{ t,y}(dx)]\\
	= \IZ (\floo yx +\floo xy ) \,\qq_t(dx) \Jf _{ t,x}(dy)
	=0,
 \end{multline*}
 where we used \eqref{eq-65} and $\floo yx +\floo xy =0.$
\end{enumerate}
\end{remarks}

Identity  \eqref{eq-65} expresses the equality of the forward and backward instantaneous  fluxes at each time $t$ and between any pair of locations $(dx,dy).$ This property which is intuitively expected when playing the movie backward, is widely used without proof in theoretical physics. Nevertheless, it appears that finding a large set of regularity assumptions on the path measure for this identity to be verified is not as easy as it seems. The aim of next Sections \ref{sec-regular-kernel} and \ref{sec-entropy} is to identify some assumptions which are more explicit.

\subsection*{Proof of Proposition \ref{res-33}}
Let us denote 
\begin{align*}
\nu(dxdy)&:=  \sigma(x,y)\qq(dx)\Jf _{ x}(dy),\\ 
\mu(dxdy)&:= \sigma(x,y) \qq(dy)\mathcal{J} _{ y}(dx).
\end{align*}
By assumption, $ \nu$ is a finite measure. Hence, its marginals are finite measures (without the boundedness of $\nu$  implied by the introduction of the function $ \sigma$, the marginals of $\qq(dx)\Jf _{ x}(dy)$ might take infinite values). In particular its $y$-marginal is $$n(dy):=\nu(\ZZ\times dy)=\If{\qq} \sigma (dy).$$
Multiplying both sides of  equation \eqref{eq-FE} by the non-vanishing function   $ \sigma,$ gives the equivalent equation
\[\nu= \mu.\] 

Suppose that \eqref{eq-FE} admits a solution $\mathcal{J}$. Then $\mu=\nu$ is a finite measure and  its
 $y$-marginal $m(dy):= \mu(\ZZ\times dy)$ is well defined. By definition of $ \mu$,   $m\ll \qq$ and  taking the $y$-marginal of  $\nu=\mu$, we obtain $n=m$  and see that \eqref{eq-70} is satisfied.

Conversely, suppose that \eqref{eq-70} holds, that is: $
 n\ll \qq.
$
Then
\begin{align*}
\nu(dxdy)
	=n(dy) \nu(dx\mid y) 
	= \frac{dn}{d\qq}(y) \qq(dy) \nu(dx\mid y),
\end{align*}
showing that
\begin{align*}
\mathcal{J} _{ y}(dx):=  \sigma(x,y) ^{ -1}\frac{dn}{d\qq}(y) \nu(dx\mid y)
\end{align*}
solves \eqref{eq-FE} uniquely. This implies the first equality in \eqref{eq-83} because $\nu(dx\mid y)$ is a probability kernel. 
\\
Let us prove the convolution expression of \eqref{eq-83}.
For all bounded measurable function  $f$ on $\ZZ$,
\begin{multline*}
\IZ f(y)\,\If \qq \sigma(dy)
	=\IZZ f( \tau^\xi(x)) \sigma(x,\tau^\xi(x))\qq(dx)\Kf_x(d\xi)\\
	=\IZZs f(y)\ \tau^\xi\pf\qq(dy)\ \sigma( \tau ^{ -\xi}(y),y)\Kf _{\tau ^{ -\xi}(y)}(d\xi).
\end{multline*}
We see that
 $\If \qq \sigma(dy)
 =\int _{ \xi\in\Zs} \sigma(y-\xi,y) \Kf _{ y-\xi}(d\xi)\ \qq^\xi (dy).$ 
With \eqref{eq-70}: $\If \qq \sigma\ll\qq,$ this implies that $\qq^\xi\ll \qq$ for almost every $\xi$  with respect to the $\xi$-marginal $\qq \sigma\Kf(\ZZ\times d\xi)$ of $\qq \sigma\Kf(dy,d\xi):= \qq(dy) \sigma(y-\xi,y)\Kf _{ y-\xi}(d\xi),$ and the second equality in \eqref{eq-83} follows. 
\hfill$\square$

\subsection*{Proof of Theorem \ref{res-37}}
It is mainly the consequence of two preliminary results:  Lemmas \ref{res-32} and  \ref{res-35} below.
\begin{itemize}
\item
Statement (a) is an immediate consequence of the local boundedness of $\Btf$ and  estimate \eqref{eq-84}. 
\item
Statement (b): $[\eqref{eq-93}\Longleftrightarrow \eqref{eq-70b}]$ is already established at Proposition \ref{res-33},  $ [\eqref{eq-93} \Longrightarrow \eqref{eq-86} ]$ is proved at Lemma \ref{res-32} whose proof uses stochastic calculus, and $[\eqref{eq-86}\Longrightarrow \eqref{eq-93} ]$ is proved at Lemma \ref{res-35} which relies essentially on the IbP formula of Theorem \ref{res-02}.
\item
Statement (c) is part of Lemma \ref{res-32}.
\end{itemize}
It  remains to state and prove Lemmas \ref{res-32} and  \ref{res-35}.

\begin{lemma} \label{res-32}
Suppose that  $\Btf$ is locally bounded, \eqref{eq-84} holds 
 and the flux equation $\FE(\qq_t,\Jf^Q_t)$ written at \eqref{eq-FE} admits a solution for almost every $t$.\\ Denote this solution by $\Jb^Q_t$ and define  $\Btb$ by \eqref{eq-67}.\\
Then, the right hand side of \eqref{eq-67} is well defined, and any $u\in C_c^2(\ZZ) $ is in $dom \LLb^Q$ with $\LLb^Q_tu$ given at \eqref{eq-85}, 
and both $\LLf^Qu$ and $\LLb^Qu$ are $\qb$-integrable.
\end{lemma}

\begin{proof}
For any $0\le r\le s\le t\le 1$ and any bounded measurable function $f$ on $\ZZ,$ we put
\begin{align}\label{eq-66}
A&:=E _{ Q^*} \Big[ f(X _{ 1-t}) \Big\{ u(X _{ 1-r})-u(X _{ 1-t})
	- \int _{ 1-t} ^{ 1-r} \Btb _{ 1-\tau}\cdot \nabla u(X_ \tau)\, d \tau\Big\} \Big] \nonumber\\
 &\ = -E_Q \left[  \left\{ u(X_t)-u(X_r) + \int _r^t \Btb_s\cdot \nabla u(X_s)\, ds \right\}f(X_t) \right]\end{align}
and 
\begin{align*}
B&:= E _{ Q^*}\Big[ f(X _{ 1-t})\Big\{ \int _{ 1-t }^{ 1-r} \big[u(x)-u(X_ \tau)-\nabla u(X_ \tau)\cdot \floo x{X_ \tau}	\big]\,\Jb _{ 1- \tau,X_ \tau}(dx) d \tau\Big\} \Big]\\
	&\ =  E_Q \left[ \left\{ \int_r^t [u(x)-u(X_s)-\nabla u(X_s)\cdot \floo x{X_s}]\,\Jb _{ s,X_s}(dx) ds\right\} f(X_t) \right]  \\
 &\ = \int _{ [r,t]\times (\ZZ)^3} [u(x)-u(y)-\nabla u(y)\cdot \floo xy ]f(z)\, \Jb _{ s,y}(dx) \qq _{ st}(dydz)ds \nonumber
\end{align*}
It will be seen during the proof that  $A$ and $B$  are well defined integrals.   
We have to prove 
\begin{align*}
A=B.
\end{align*}
Because $f(X_t)$ depends on the future of the remaining terms of the integrands in $E_Q$, we are in a bad shape to attack this problem with martingale techniques. In fact, we shall rely on Itô's formula 
\begin{align}\label{eq-68}
\begin{split}
u(X_t)-u(X_r)
	= \sum _{ s\in[r,t]}[u(X_s)-u(X _{ s^-})&-\nabla u(X _{ s^-}) \cdot \floo{X_s}{X _{ s^-}}]\\
	&+ \int_r^t \Btf_s(X_s)\cdot\nabla u(X_s)\,ds,
		\qquad Q\ae,
\end{split}
\end{align}
which is an almost sure identity. If the sample paths have bounded variations, the series $\sum _{ s\in[r,t]} \cdots$ is defined in the usual sense. In the general case, it is a stochastic integral whose compensator with respect to $Q$  is 
\begin{align*}
\int_r^t [u(y)-u(X _{ s^-})-\nabla u(X _{ s^-})\cdot \floo y{X _{ s^-}}]\, \Jf _{ s, X _{ s^-}}(dy)ds.
\end{align*}
With  \eqref{eq-65} one obtains
\begin{align*}
\Jb _{ s,y}(dx) \qq _{ st}(dydz)
	=\Jb _{ s,y}(dx) \qq _{ s}(dy)\qq_t(dz\mid X_s=y)
	=\qq_s(dx)\Jf _{ s,x}(dy)\qq_t(dz\mid X_s=y),
\end{align*}
which transforms $\Jb _{ s,y}(dx) \qq _{ st}(dydz)ds$ whose meaning  is obscure, into the meaningful  (at least when $\Jf(\ZZ)< \infty$) expression
\begin{align}\label{eq-69}
\begin{split}
\qq_s(dx)&\Jf _{ s,x}(dy)\qq_t(dz\mid X_s=y)ds \\
	&\simeq Q \big( X _{ s^-}\in dx,\  \textrm{there is a jump $x\to dy$ during }(s,s+ds],\  X_t\in dz \big) .
\end{split}
\end{align}
This leads us to
\begin{align*}
B= \int _{ [r,t]\times (\ZZ)^3} [u(x)-u(y)-\nabla u(y)\cdot \floo xy ]f(z)\, \qq_s(dx)\Jf _{ s,x}(dy)\qq_t(dz\mid X_s=y)ds
\end{align*}
and proves that $B$  is a well defined integral under our integrability assumption  \eqref{eq-84}. \\
Let us force the appearance of $\nabla u(x)$ in place of $\nabla u(y).$
With  \eqref{eq-65} again, we see that
\begin{multline}\label{eq-74}
\int _{[r,t]\times \ZZZ} \{\nabla u(y)-\nabla u(x)\} \cdot \floo yx 
\ \qq_s(dx)\Jf _{ s,x}(dy)ds\\
	=- \int _{ [r,t]\times\ZZ} \nabla u(y)\cdot c_s(y)\, \qq_s(dy)ds
\end{multline}
where we set 
\begin{align*}
c_s(y):= \IZ \floo xy  \,[\Jf _{ s,y}+\Jb _{ s,y}](dx).
\end{align*}
Note in passing that  $c$ is well defined and integrable with respect to $\qb,$ because  the integral  on the left hand side is finite.
It follows that 
\begin{align*}
B &= -\int _{ [r,t]\times (\ZZ)^3} [u(y)-u(x)-\nabla u(x)\cdot \floo yx ]f(z)\, \qq_s(dx)\Jf _{ s,x}(dy)\qq_t(dz\mid X_s=y)ds\\
	&\hskip 3cm -\int _{ [r,t]\times (\ZZ)^3}  \nabla u(y) \cdot c_s(y)\, f(z)
\ \qq_s(dy)\qq_t(dz\mid X_s=y)ds\\
 &= -E_Q \Big[ \Big\{ \sum _{ s\in [r,t]}[u(X_s)-u(X _{ s^-})-\nabla u(X _{ s^-})\cdot \floo{X_s }{X _{ s^-}}\Big\}  f(X_t)\Big] \\
 &\hskip 8cm -E_Q\Big[f(X_t) \int_r^t\nabla u(X_s)\cdot c_s(X_s)\,ds\Big] ,
\end{align*}
where the main idea  for  last identity is   \eqref{eq-69}, but it is valid even when the jump frequency is infinite.
With Itô's formula \eqref{eq-68},
 we arrive at
\begin{align*}
B= -E_Q \left[  \left\{ u(X_t)-u(X_r) 
	+ \int_r^t (c_s-\Btf_s)(X_s)\cdot \nabla u(X_s)\,ds\right\} f(X_t)\right] .
\end{align*}
Going back to the expression \eqref{eq-66} of $A$,  we see that the desired identity $A=B$ is realized once $\Btf+\Btb=c,$ which is \eqref{eq-67}. \end{proof}

Recall that  Theorem \ref{res-02} states that  any function $u$ in $\UQ$ such that 
$u\in\dom \LLb^Q$,   $\LLb^Qu\in L^1(\qb)$,  verifies 
the IbP formula \eqref{eq-24}  for almost every $t$:
\begin{align*}
\IZ [(\LLf^Q_tu+\LLb^Q_tu)v + \Gaf^Q_t(u,v)]\, d\qq_t=0,
\quad \forall v\in\UQ.
\end{align*}

\begin{lemma}\label{res-35}
Assume that $Q\in\MO$ is Markov and $\UQ$ is measure determining.
\\
Then, for any $u\in\UQ\cap \dom\LLb^Q$ such that $\LLb^Qu\in L^1(\qb),$ the backward generator writes as
\begin{align*}
\LLb^Q_tu(x)
	=\Btb(x)\scal\nabla u(x)
		+\IZ [u(y)-u(x)-\nabla u(x)\scal \floo yx]\,\Jb _{ t,x}(dy),
\end{align*}
where $\Jb$ solves equation \eqref{eq-65} and $\Btb$ is defined by \eqref{eq-67}.
\end{lemma}

\begin{proof}
Let us start the calculations in the simplest case where there is no drift and the sample paths have bounded variations:
\begin{align*}
\LLf u(x)=\IZ [u(y)-u(x)]\, \Jf_x(dy),
\end{align*}
where we dropped the time subscript $t$ for simplicity. By Theorem \ref{res-02} the IbP formula holds:
\begin{align}\label{eq-73}
\IZ& \LLb u\, v\, d\qq
	=-\IZ \LLf u\, v\,d\qq -\IZ \Gaf (u,v)\,d\qq \\
	&=-\IZ[u(y)-u(x)] v(x)\,\qq(dx)\Jf_x(dy)-\IZ[u(y)-u(x)][v(y)- v(x)]\,\qq(dx)\Jf_x(dy)\nonumber\\
	&=-\IZ[u(y)-u(x)] v(y)\,\qq(dx)\Jf_x(dy),\nonumber
\end{align}
for all $u,v\in\UQ.$ Hence,
\begin{align}\label{eq-72}
\LLb u(x)\qq(dx)
	&=\IZ[u(y)-u(x)] \,\qq(dy)\Jf_y(dx)\\
	&= \left(\IZ [u(y)-u(x)]\kappa(dy|x)\right) \ \rr(dx)\nonumber
\end{align}
where we set $\kappa(dxdy):=\qq(dy)\Jf_y(dx)=\rr(dx)\kappa(dy| x)$ with $\kappa(dy| x)$ a probability kernel.
Setting $f^u(x):=\LLb u(x),$ $g^u(x):=\IZ [u(y)-u(x)]\, \kappa(dy|x)$, this identity writes as
\begin{align*}
f^u(x)\qq(dx)=g^u(x)\rr(dx).
\end{align*}
Because $\1 _{ g^u\neq 0}\, \rr\ll g^u\rr=f^u\qq\ll\qq$, we can write 
$f^u\qq=g^u(\1 _{ g^u\neq 0}\rr)=g^u  \frac{d(\1 _{ g^u\neq 0}\rr)}{d\qq}\qq,$ showing that
\begin{align*}
f^u=  \frac{d(\1 _{ g^u\neq 0}\rr)}{d\qq}\ g^u,\quad \qq\ae
\end{align*}
Beware,   we look for a formula $f^u= \alpha g^u$ where the function $ \alpha$ does not depend on $u$. However, for a countable subclass $\widetilde\UU$ of $\UQ$, the set $S:=\cup _{ u\in\widetilde\UU}\{g^u\neq 0\}$ is measurable,  $  \1 _S\, \rr\ll  \qq,$ and 
\begin{align*}
f^u= \alpha g^u,\quad \qq\ae,\quad \forall u\in\widetilde\UU,
\quad \textrm{where }\alpha= \frac{d(\1 _S\rr)}{d\qq}.
\end{align*}
This proves that for any $u\in\widetilde \UU$, 
$\LLb u(x)\qq(dx)=\IZ \IZ [u(y)-u(x)]\, \alpha(x)\kappa(dy|x)\qq(dx),$ which gives
\begin{align*}
\LLb u(x)=\IZ[u(y)-u(x)]\,\Jb_x(dy),\quad\forall x\ \qq\ae,
\quad \forall u\in\widetilde\UU,
\end{align*}
with $\Jb_x(dy)= \alpha(x)\kappa(dy|x).$ The martingale problem associated to $(\LLb, \UQ)$ is completely specified by its restriction to a large enough countable subclass $\widetilde\UU$ of $\UQ,$ because the $ \sigma$-field on $\OO$ is countably generated. It follows that   the above expression of $\LLb u$ extends to any $u$ in $\UQ.$ With \eqref{eq-72} we arrive at
\begin{multline*}
\IZZ v\LLb u\,d\qq=\IZZ[u(y)-u(x)] v(x)\,\qq(dy)\Jf_y(dx)\\
=\IZZ[u(y)-u(x)]v(x)\,\qq(dx)\Jb_x(dy).
\end{multline*}
This leads us to the flux identity \eqref{eq-65} because the collection of all functions $(x,y)\mapsto [u(y)-u(x)]v(x)$ when $u$ and $v$ describe $\UQ$ is measure-determining on $\ZZZ$ outside the diagonal.

We now look at the case where
\begin{align*}
\LLf u(x)=\Bf(x)\scal \nabla u(x)+ \IZ [u(y)-u(x)]\,\Jf_x(dy)
\end{align*}
where $\Bf$ is  locally bounded and the sample paths have bounded variations. The only difference with previous computations is the addition of $\Bf(x)\scal \nabla u(x)$. Going back to \eqref{eq-73}, we see that we have to replace $\LLb u $ by $\LLb u +\Bf\scal \nabla u,$ leading us to the same flux identity \eqref{eq-65}, and to $\Bb=-\Bf.$ 

In the general case where $\LLf$ is given at \eqref{eq-752},  one completes the proof proceeding as in previous calculations and reasoning as in the proof of Lemma \ref{res-32} at \eqref{eq-74}. 
\end{proof}

\section{Regular jump kernel} 
\label{sec-regular-kernel}

\subsection*{Time reversal  without   IbP}

Next result is a time reversal formula which does not rely on the IbP formula.

\begin{proposition}\label{res-34}
Suppose that $Q\in\MO$ satisfies \eqref{eq-84},  $\UQ= C_c^2(\ZZ),$ and   $\Btf_t$,  $\Kf^Q _{ t}$ are in $C^1_x$ for all $t>0$. 
\\
Suppose also that 
 $\qq_t$ is absolutely continuous for all $t>0$ and the flow of densities  $(t,x)\mapsto q_t(x):=d\qq_t/dx$ is in $C ^{ 1,1} _{ t,x}.$
\\
Then, the flux equation $\FE(\qq_t,\Jf_t^Q)$ written at \eqref{eq-FE} admits a solution $\Jb_t^Q$  for all $t>0,$  and the time reversal formula \eqref{eq-85}-\eqref{eq-65}-\eqref{eq-67} is valid.
\end{proposition}

\begin{proof}
By Proposition \ref{res-33} and Lemma \ref{res-32} which do not rely on the IbP formula, it is sufficient  to show that \eqref{eq-70b} holds for all $t>0:$   we have to prove
\begin{align*}
\qq_t(A)=0\implies \If{t,\qq_t} \sigma(A)=0,
\end{align*}
for all $t>0.$
The  evolution of $\qq_t$ is governed by the weak equation
\begin{multline*}
\frac{d}{dt} \langle u, \qq_t\rangle 
	=  \IZ \Btf (t,x)\scal\nabla u(x) q_t(x)\, dx  \\
	+ \IZZs [u(y)-u(x)-\nabla u(x)\cdot \floo yx ]q_t(x)\, dx\Jf_x(dy)
\end{multline*}
for any $u\in C^2_c(\ZZ).$ Since $q_t, \Btf$ and $\Kf$ are assumed to be differentiable in space, integrating by parts we obtain
\begin{multline*}
\frac{d}{dt} \langle u, \qq_t\rangle 
	=- \IZ u(x) \ \nabla\scal (q_t\Btf _t)(x)\, dx\\
	 +\IZ u(x)\, ``[\nu_t(\ZZ\times dx)- \nu_t(dx\times \ZZ)+\nabla\scal (q_t\beta^\delta)(x)\, dx]''
\end{multline*}
where $ \nu_t(dxdy):=\qq_t(dx)\Jf _{ t,x}(dy)$. It is understood that under the assumption 
$	%\begin{align*}
\IZs (1\wedge |\xi|)\, \Kf(d\xi)< \infty,
$	%\end{align*}
$\beta^\delta:=\IZs \xd\, \Kf(d \xi)$ and   
\[
\nu(\ZZ\times dx)- \nu(dx\times \ZZ)
	=\nu( [dx]^c\times dx )- \nu( dx\times [dx]^c)
\]
with
$[dx]^c:=\{\ZZ\setminus dx\}$ (remark that the diagonal is not charged, as desired). Under the general hypothesis $	%\begin{align*}
\IZs (1\wedge |\xi|^2)\, \Kf(d\xi)< \infty,
$
this expression must be compensated by  the contribution $-\nabla\scal (q_t \beta^\delta)$ of the small jumps appearing in the ill-defined integral $\beta^\delta=\IZs \xd\, \Kf(d \xi).$ In this case, none of each separate terms of  $[\nu_t(\ZZ\times dx)- \nu_t(dx\times \ZZ)+\nabla\scal (q_t\beta^\delta)(x)\, dx]$ is meaningful, contrary to the whole expression.

  This equation extends to any    integrable measurable test function $u$ with a  bounded support. In particular, with $u=\1_A$ the indicator of a bounded measurable subset $A$ satisfying $\qq_t(A)=0$, we obtain
\begin{align*}
\frac{d}{dt}\qq_t(A)
	= n_t (A),
\end{align*}
with $ n_t (A):= \nu_t(\ZZ\times A)\ge 0.$
To see this, remark that  for $\qq_t$-almost every $x$ in $A$, we have $q_t(x)=0$ and $\nabla q_t(x)=0$ since $q_t$ is assumed to be differentiable and $q_t(x)=0$ is a minimal value. Hence,  the divergence integrals vanish: $\int_A \nabla\cdot (q_t v)(x)\,dx= \int_A \nabla v\, d\qq_t+\int_A \nabla q_t\cdot v(x)\, dx=0.$ On the other hand, $\nu_t(A\times \ZZ)=0$ because $\nu_t( \sbt \times \ZZ)\ll\qq_t.$ Hence the only remaining term in the right hand side is $\nu_t(\ZZ\times A)=:n_t(A).$

Supposing ad absurdum that $ n_t (A)>0$ implies that $\qq_s(A)<0$ for some $0\le s<t,$ a contradiction. Therefore, $ n_t (A)=0,$ which in turns implies that $\If{t,\qq_t} \sigma(A)=0$ for any positive $ \sigma.$
\end{proof}

%\subsection*{Time reversal formulas based on the  IbP formula}

The hypotheses of Proposition \ref{res-34} are rather restrictive. In particular, any Poisson process starting from a Dirac  mass is ruled out by the requirement that the time marginals are absolutely continuous with respect to Lebesgue measure. 

In contrast  Theorem \ref{res-36} below    offers us  more handy sets of assumptions for the time reversal formula.
 Unlike previous Proposition \ref{res-34}, its proof does not rely on Lemma \ref{res-32} and the resolution of equation \eqref{eq-65}, but on  the IbP formula. Moreover, the existence of a solution to \eqref{eq-65} is obtained as a corollary.

\subsection*{A time reversal formula based on Theorem \ref{res-02}-(b)}

Let us prove  a time reversal formula based on the IbP formula of Theorem \ref{res-02}-(b) when the forward jump kernel is regular enough for the  carré du champ $\Gamma^Q$  to verify  \eqref{eq-39b} and \eqref{eq-39c}. 

\begin{theorem}\label{res-36}
Suppose that $Q\in\MO$  satisfies the Hypotheses \ref{ass-08}-a-b-c.\\
Then, any  $u$ in $ C_c^2(\ZZ)$ is in the domains of $\LLf^Q$ and $\LLb^Q$, the equation \eqref{eq-65} admits a solution $\Jb^Q$ and $\LLb^Q$ is given by \eqref{eq-85} with $\Btb$  defined by \eqref{eq-67}.
\end{theorem}

\begin{remarks}\ 
\begin{enumerate}[(a)]

\item
Hypothesis \ref{ass-08}-d is not required here. It will be necessary for the proof of Theorem \ref{res-43}.

\item
The  assumptions  of Theorem \ref{res-36} are less restrictive than those of Proposition \ref{res-34}. This is obvious when the sample paths have finite variation. Otherwise, when the sample paths have infinite variation,  the assumed $t$-differentiability of $q$ in Proposition \ref{res-34} is replaced by some $x$-differentiability of the jump kernel which is easier to verify.

\item
Instead of the full statement  \eqref{eq-88} of Hypothesis (c2) which requires that $q$ is $C _{ t,x} ^{ 0,1}$, the proofs of Theorem \ref{res-36} and its preliminary Lemma \ref{res-41}  only rely on:  $q_t$ is $C^1_x$ for all $t$. 
\end{enumerate}

\end{remarks}

\begin{proof}[Proof of Theorem \ref{res-36}] \label{pf-res36}
This result is a direct corollary of  Thm.\ \ref{res-02} and Lemma \ref{res-35}.  All we have to do is to make sure that the hypotheses of Thm.\ \ref{res-02}-(b) are verified.  

It is easy to see that under the general assumptions of the proposition:  \eqref{eq-32}, \eqref{eq-94} and    $[\eqref{eq-76}\ \mathrm{or}\ \eqref{eq-110}],$  we have  \eqref{eq-84}  and  $\UQ=\UU,$ that is: $ \LLf^Qu \in L^1(\qb)$ and  $ \Gaf^Q (u) \in L^1(\qb) $ for any $u\in \UU:=C_c^2(\ZZ),$ recall \eqref{eq-42}.    
Note that $C_c^2(\ZZ)$ determines the weak convergence of Borel measures on $\ZZ,$ as desired. It is also clear that \eqref{eq-39b} holds, i.e.\ $(t,x)\mapsto \Gaf^Q_t(u,v)(x)$ is continuous, under the regularity assumptions \eqref{eq-771} or \eqref{eq-772}.

It remains to verify \eqref{eq-39c}, that is: For any $u\in C_c^2(\ZZ),$ the linear form
\begin{multline}\label{eq-120}
w\in C_c ^{ 0,2}(\iZ)\\ \mapsto \ell(w):= \int _{ \ii\times\ZZ\times\Zs} [w_t(x+\xi)-w_t(x)][u(x+\xi)-u(x)]\qq_t(dx) \Kf^Q  _{ t,x}(d\xi)dt
\end{multline}
 is a finite signed measure on $\iZ$. Denoting the measure
 \begin{align*}
  \mu(dtdxd\xi):=[u(x+\xi)-u(x)]\qq_t(dx) \Kf^Q  _{ t,x}(d\xi)dt
 \end{align*}
and  the mapping $$ \tau(t,x,\xi):=(t,x+\xi, -\xi)$$ on $ \ii\times\ZZ\times\Zs,$ we see that
\begin{align*}
\ell(w)=\int _{ \ii\times\ZZ\times\Zs} w(t,x)\  [\tau\pf \mu- \mu](dtdxd\xi)
	=\IiZ w(t,x)\  [\tau\pf \mu- \mu](dtdx\times \Zs),
\end{align*}
where these identities are formal. Indeed, when $\Kf^Q (\Zs)= \infty,$  the  term $[\tau\pf \mu- \mu](dtdx\times \Zs)$ might not even be defined as a measure. To complete the proof of the proposition, we have to show that under our assumptions, 
$$
[\tau\pf \mu- \mu](dtdx\times \Zs)
$$ 
is a finite measure on $\iZ.$

Each assumption \eqref{eq-76} or \eqref{eq-110}  implies that $ \mu(\iZ\times B_ 1^c)< \infty$ and also that  $ \tau\pf\mu(\iZ\times B_ 1^c)= \mu(\iZ\times B_ 1^c)< \infty$  (remark that $-B^c_ 1=B^c_ 1$).

\subsubsection*{Bounded variation case}

Under the assumptions (i) corresponding to the bounded variation case, we see that  $ \mu(\iZ\times B_ 1)< \infty$ because  $u(x+\xi)-u(x)$ is close to $\nabla u(x)\scal \xi$ for small jumps $\xi$ and that $\nabla u$ is bounded. Again, we have $\tau\pf \mu(\iZ\times B_ 1)= \mu(\iZ\times B_ 1)< \infty.$

\subsubsection*{Unbounded variation case} 
Under the assumptions (c2) it is proved  at Lemma \ref{res-41} below that $\ell$ is a bounded measure.
\end{proof}

The proof of the remaining  Lemma \ref{res-41} relies upon the preliminary Lemma \ref{res-40} below. 

By hypothesis, 
 the jump
\begin{align*}
\phi _{ t, \alpha} (x)= T _{ t, \alpha}(x)-x=  \nabla\theta _{ t, \alpha}(x) -x= \nabla \psi _{ t, \alpha}(x),
\end{align*}
writes as a displacement from $x$ to $T _{ t, \alpha}(x)$ with $T _{ t, \alpha}=\nabla \theta _{ t, \alpha}$
for some function  $$ \theta _{ t, \alpha}(x)= |x|^2/2+ \psi _{ t, \alpha}(x), \qquad x\in\ZZ,$$ which is strictly convex and differentiable, by assumption \eqref{eq-100b}. This is the well-known framework of quadratic optimal transport where $x\mapsto T _{ t, \alpha}$  is a Brenier mapping, see \cite{Vill09}.
\\
Let us fix $t, \alpha$ for a while and drop the indices $t, \alpha.$ Under the above assumption the Brenier mapping $T$ is invertible and
\[
y=T(x)= \nabla \theta(x) \iff x=T ^{ -1}(y)=\nabla \theta^*(y)
\]
where $ \theta^*$ is the convex conjugate of $ \theta.$ Therefore
\begin{align}\label{eq-98}
\begin{split}
&y-x=\phi(x) =\nabla  \psi(x),
\qquad \psi(x)= \theta(x)-|x|^2/2,\\
&x-y=\hat\phi(y) =\nabla \hat \psi(y),\qquad \hat\psi(y)= \theta^*(y)-|y|^2/2.
\end{split}
\end{align} 
We interpret $\phi(x)$ as the forward jump from $x$ to $y$, and $\hat\phi(y)$ as the backward jump from $y$ to $x$.
The regime we investigate is $\phi(x) $ close to zero. We assume that the function $\psi$ is $C^2$, and satisfies
\begin{align*}
\nabla \psi(x)\neq 0,\quad \forall x\in\ZZ,
\quad \textrm{and} \quad
\nabla^2 \psi\ge c''\Id\   \textrm{ for some }\   c''\in\RR.
\end{align*}
 Let us give a name to the bounds on the derivatives of $\psi$: 
\begin{align*}
\sup _{x\in\ZZ}|\nabla \psi(x)|=: C',\qquad
\sup _{B _{ \rho+C'}}|\nabla^2 \psi|=: C''_ \rho,
\end{align*}
where  $B_ r:= \left\{ z\in\ZZ; |z| \le r\right\}$ denotes the ball  centered at zero with radius $r\ge 0$.

\begin{lemma}\label{res-40}
Assume that  $c''>-1$ and  $C'< \infty$.\\ Then, for any $x\in\ZZ,$ 
\begin{align}\label{eq-102}
\hat\phi(x)&=- \phi(x+z_x),\\
\nabla\hat\phi(x)&=\nabla \phi(x+z_x). \label{eq-106}
\end{align}
where $z_x$ is the unique solution of 
\begin{align}\label{eq-101}
z_x=-\phi(x+z_x).
\end{align}
Moreover, for any $\rho>0$  and all $x\in B_  \rho,$
\begin{align}\label{eq-105}
0<1+  c''
	\le \frac{|\phi(x)|}{|\hat\phi(x)|}
	\le 1+ C''_\rho .
\end{align}
Suppose that in addition, for all $ \rho\ge 0,$ there exists $c_ \rho< \infty$ such that
\begin{align}\label{eq-100db}
|\nabla\phi(x)|\le c_ \rho|\phi(x)|,
\qquad \forall x\in B _{ \rho+C'}.
\end{align}
Then, for any $ \rho>0$ and all $x\in B_ \rho,$
\begin{align}\label{eq-103}
|\phi(x)+\hat\phi(x)|
	 \le  \frac{c_ \rho(1+c''+c_ \rho C') }{(1+c'')^2}\   |\phi(x)|^2.
\end{align}
\end{lemma}
The hypotheses 
  $c''>-1$,  $C'< \infty$ and \eqref{eq-100db} correspond to the hypotheses \eqref{eq-100b}, \eqref{eq-100c} and \eqref{eq-100d}.

\begin{proof}
 With \eqref{eq-98}
\begin{align*}
\hat\psi(x)
	&=\theta^*(x)-|x|^2/2
	=\sup _{ y} \left\{ y\scal x - \theta (y) \right\} -|x|^2/2\\
	&=\sup _{ y} \left\{ y\scal x - \psi (y)-|y|^2/2-|x|^2/2\right\} 
	=-\inf _{ y} \left\{ |x- y|^2/2 + \psi (y)\right\} \\
	&=-\inf _{ z} \left\{ | z|^2/2 + \psi (x+ z)\right\}
%	=- \psi(y)-\inf _{ z} \left\{ | z|^2/2 +  \nabla \psi(y)\scal z+ r(y,z)\right\}.
\\
	&= - \psi(x)- \inf _{ z} \left\{ | z|^2/2 + [\psi (x+ z)-\psi(x)]\right\}.
\end{align*}
Since  $\nabla^2\psi\ge c''\Id>-\Id$,  the function $z\mapsto | z|^2/2 + [\psi (x+ z)-\psi(x)]$ is strictly convex and it achieves its unique minimum at $z_x$, solution of \eqref{eq-101}: $z_x=-  \nabla \psi(x+z_x)$.
This { implies that}
\begin{align}\label{eq-99}
 |z_x|\le C' ,\qquad \forall x\in \ZZ
\end{align}
with $C'< \infty$ by hypothesis, and 
\begin{align*}
\psi(x)+\hat\psi(x)
=
-\inf_z \cdots
	=- |\nabla \psi(x+z_x)|^2/2
		-   [\psi (x+ z_x)-\psi(x)].
\end{align*}
Differentiating once more
\begin{align*}
\phi(x)&+\hat\phi(x)\\
	&=\nabla\psi(x)+\nabla\hat\psi(x)\\
	&= -(\Id+ z'_x)^T\nabla^2\psi(x+z_x)\nabla \psi(x+z_x)
	- \Big[(\Id+z'_x)^T \nabla \psi(x+z_x)-\nabla\psi(x) \Big]\\
	&=- \nabla^2\psi(x+z_x)\nabla \psi(x+z_x)
	 - \nabla\phi(x+\bar z_x)\, z_x\\
	&\hskip 3cm -  z_x ^{ \prime T}\, \nabla\psi(x+z_x) 
	-   z_x ^{ \prime T}\,\nabla^2\psi(x+z_x)\nabla \psi(x+z_x)\\
	&=- \nabla^2\psi(x+z_x)\nabla \psi(x+z_x)
	+\nabla^2\psi(x+\bar z_x)\, \nabla\psi(x+z_x)\\ 
	&\hskip 3cm - z_x ^{ \prime T}\, \nabla\psi(x+z_x)
	-   z_x ^{ \prime T}\,\nabla^2\psi(x+z_x)\nabla \psi(x+z_x)\\
	&=\nabla\phi(x+\bar z_x)\, \phi(x+z_x)\\
	&\hskip 2cm - \nabla\phi(x+z_x)\phi(x+z_x)
	-  z_x ^{ \prime T}\,(\Id+  \nabla\phi(x+z_x)) \phi(x+z_x)
\end{align*}
for some $\bar z_x\in[0,z_x]$ because $\psi$ is assumed to be $C^2.$  The derivative $z'_x$ of $x\mapsto z_x$ exists by  local inversion,  because $c''>-1$  implies that  $\Id+ \nabla^2\psi=\Id+\nabla\phi$ is invertible,  and
\begin{align}\label{eq-107}
z'_x
	=- (\Id+ \nabla\phi(x+z_x)) ^{ -1}\,\nabla\phi(x+z_x).
\end{align}
This gives
\begin{multline*}
-  z_x ^{ \prime T}\,(\Id+  \nabla\phi(x+z_x)) \\
	= \Big[\nabla\phi\, \{\Id+ \nabla\phi \} ^{ -1}
	\,\{\Id+  \nabla\phi\}\, \phi \Big](x+z_x)
	=  \nabla\phi(x+z_x)\, \phi(x+z_x).
\end{multline*}
Hence
\begin{align}\label{eq-102b}
\phi(x)+\hat\phi(x) 
	=\nabla\phi(x+\bar z_x)\, \phi(x+z_x).
\end{align}
On the other hand,
\begin{align*}
\phi(x+z_x)-\phi(x)
	= \nabla\phi(x+\bar z_x)\,z_x
	=-\nabla\phi(x+\bar z_x)\, \phi(x+z_x)
\end{align*}
with the same $\bar z_x\in[0,z_x]$ as above. Comparing with \eqref{eq-102b} proves \eqref{eq-102}, which in turns implies \eqref{eq-106} because   
\begin{align*}
\nabla\hat\phi(x)
	&= -\nabla\phi(x+z_x)\,(\Id+z'_x)\\
	&\overset{\eqref{eq-107}}=-\nabla\phi(x+z_x)\,\Big(\Id-\big(\Id+\nabla\phi(x+z_x)\big)\Big)^{ -1}\,\nabla\phi(x+z_x)
	=\nabla\phi(x+z_x).
\end{align*}
To prove \eqref{eq-105}, we see with \eqref{eq-102} and  \eqref{eq-102b}  that 
\begin{align}\label{eq-123}
\hat\phi(x)+\phi(x)=-\nabla \phi(x+\bar z_x)\hat\phi(x)\iff 
-\phi(x)=(\Id +\nabla \phi(x+\bar z_x))\hat\phi(x).
\end{align}
It follows that for all $x\in B_ \rho,$
\begin{align*}
1+c''\le
\inf _{ x'\in B _{ \rho+C'}} \|\big(\Id+\nabla\phi(x')\big) ^{ -1}\| ^{ -1}
	\le
\frac{|\phi(x)|}{|\hat\phi(x)|}
	\le \sup _{ x'\in B _{ \rho+C'}} \|\Id+\nabla\phi(x')\|
	\le 1+ C''_ \rho.
\end{align*}
It remains to prove \eqref{eq-103}. Again, the starting point is \eqref{eq-123}. On the other hand, using our assumption \eqref{eq-100db}
\begin{multline*}
 |\nabla \phi(x+\bar z_x)|
 	\overset{ \eqref{eq-100db}}\le c_ \rho  | \phi(x+\bar z_x)|
	= c_ \rho  | \phi(x)+\nabla\phi(x+\tilde  z_x)\scal \bar z_x|
	\overset{ \eqref{eq-100db}}\le c_ \rho|\phi(x)|+ c_ \rho^2 |\phi(x+\tilde z_x)||\bar z_x|\\
	\le c_ \rho|\phi(x)|+ c_ \rho^2 C' | z_x|
	\overset{ \eqref{eq-102}, \eqref{eq-101}}= c_ \rho|\phi(x)|+ c_ \rho^2 C' | \hat\phi(x)|
	\overset{ \eqref{eq-105}}\le \big(c_ \rho +c_ \rho^2C' (1+c'') ^{ -1}\big) |\phi(x)|
\end{multline*}
and finally
\begin{multline*}
|\hat\phi(x)+\phi(x)|
	\overset{ \eqref{eq-123}}\le |\nabla \phi(x+\bar z_x)| |\hat\phi(x)|
	\overset{ \eqref{eq-105}}\le |\nabla \phi(x+\bar z_x)|(1+c'') ^{ -1} |\phi(x)|\\
	\le \big(c_ \rho +c_ \rho^2C' (1+c'') ^{ -1}\big) (1+c'') ^{ -1} |\phi(x)|^2,
\end{multline*}
which is \eqref{eq-103}.
\end{proof}

\begin{lemma}\label{res-41}
Under the Hypotheses \ref{ass-08}-a-b-c2 in the unbounded variations case where \eqref{eq-97} holds,  for any $u\in C_c^2(\ZZ)$ the linear form $\ell$ defined at \eqref{eq-120}
 is a finite signed measure on $\iZ$. 
\end{lemma}

\begin{proof}
The large jump contribution of \eqref{eq-120} integrates by parts directly. We focus on the small jump contribution.
Fix $t$ for a while and drop it as an index. We want to integrate by parts the integral
\begin{multline*}
\IZZs \cchi_ \delta(\xi) [ w(x+ \xi) - w(x)] [u(x+\xi)-u(x)] \, \qq(dx)  K_x(d \xi)\\
=\IZZs [ w(x+ \xi)-w(x)] A(x, \xi) \,\hat K_x(d \xi) dx
\end{multline*}
where we denote  for a better readability
\begin{align*}
\hat K_x(d\xi)
	&:= [\phi_x]\pf \Lambda(\xi) \\
A(x, \xi)
	&:=a(x,\xi)  k(x,\xi)q(x),\quad 
a(x,\xi)
	:=\cchi_\delta(\xi)\  [u(x+\xi)-u(x)]
\end{align*}
with   $\cchi_ \delta$ is a smooth version of the indicator of $B_ \delta$, that is:     $\cchi_ \delta$ is $C^1$ and it satisfies: $0\le \cchi_\delta\le 1$ ,  $\cchi_\delta(\xi)= \left\{ \begin{array}{ll}
1,\   & \textrm{if }|\xi|\le \delta/2,\\
0, & \textrm{if }|\xi|\ge \delta,
\end{array}\right. $    and $\sup |\nabla \cchi_\delta|\le 4/\delta$. 
We take 
\begin{align*}
0< \delta= \delta_o( \rho_u)\le 1,\quad \textrm{(recall \eqref{eq-100})}\qquad
 \textrm{with }\quad  \rho_u:= \max _{ x\in\supp u}|x|+1< \infty,
\end{align*}
to take fully account of the small jump contribution in \eqref{eq-120}.
\\
Since
\begin{align*}
\IZZs w(x) A(x, \xi) \, K_x(d \xi)dx
	=\IZA w(x) A(x, \fa(x))\, \Lambda(d \alpha)dx
\end{align*}
and
\begin{align*}
\IZZs w(x+ \xi) A(x, \xi) \, K_x(d \xi)dx
= \IZA  w(x+ \phi_ \alpha( x)) A(x, \phi_ \alpha( x))\, \Lambda(d \alpha)dx\\
=\IZA w(y) A\big(y+\faa(y), - \faa(y)\big)\ \big| \det(\Id+\nabla\faa)\big|(y)\, \Lambda(d \alpha)dy,
\end{align*}
we obtain the integration by parts formula
\begin{align}\label{eq-108}
\begin{split}
& \lambda(w):=\IZZs  [w(x+ \xi)-w(x)] A(x, \xi) \, \hat K_x(d \xi)dx\\
	&\ \ =\IZA w(x)\Big\{ A\big(x+\faa (x),-\faa (x)\big)\, (1+ \zeta_ \alpha(x)) -A\big(x, \fa(x)\big)\Big\}\, \Lambda(d \alpha)dx
\end{split}
\end{align}
where we set 
\begin{align}\label{eq-111}
| \det(\Id+\nabla\faa)|=:1+ \zeta_ \alpha.
\end{align} 
 We  need the first order expansion
\begin{align*}
A\big(x+\faa (x),-\faa (x)\big)\, (1+ \zeta_ \alpha(x)) -A\big(x, \fa(x)\big)
	=B_ \alpha(x)+B^x_ \alpha(x)+B^\xi_ \alpha(x)
\end{align*}
where
\begin{align*}
B_ \alpha(x)&:=A\big(x+\faa (x),-\faa (x)\big)\, \zeta_ \alpha(x)\\
B^x_ \alpha(x)	&:= \nabla_x A\big(x, s _{ \alpha,x}\faa(x)\big)\scal \faa(x)\\
B^\xi_ \alpha(x)		&:= -\nabla_ \xi A\big(x,(1+t _{  \alpha,x}) \fa(x)+t _{ \alpha,x}\faa(x)\big)\scal \big( \fa(x)+\faa(x)\big)
\end{align*}
for some $0\le s _{ \alpha,x} t _{ \alpha,x}\le 1,$ which is valid because $A$ is $C ^{ 1,1} _{ x, \xi}$.  We  have also
\begin{align*}
A(x,\xi) 
	&= a(x,\xi)k(x,\xi) q(x) ,\\
\nabla_x A(x, \xi)
	&= aq(x,\xi) \nabla_x k(x,\xi)
		+kq(x,\xi) \, \nabla_x a(x, \xi)+ ak(x,\xi) \nabla  q (x) ,\\
\nabla_ \xi A(x, \xi)
	&=  aq(x,\xi) \nabla_\xi k(x,\xi)
		+kq(x,\xi) \, \nabla_\xi a(x, \xi),
\end{align*}
with
\begin{align*}
|a(x,\xi)|
	&\le \sup|\nabla u|\, \1 _{B _{ \rho_u}}(x)\ \1 _{\{ |\xi|\le \delta\}}\ |\xi|,\\
|\nabla_x a(x,\xi) |&=\cchi_\delta(\xi)|\nabla u(x+\xi)-\nabla u(x)|
	\le \sup |\nabla^2 u|\ \1 _{B _{ \rho_u}}(x)\ \1 _{\{ |\xi|\le \delta\}}\ |\xi|,\\
|\nabla_\xi a(x,\xi) |
	&=\Big|\cchi_\delta(\xi)\nabla u(x+\xi)
		+  [u(x+\xi)-u(x)]\ \nabla \cchi_\delta(\xi)\Big|\\
	&\hskip 6cm \le 5 \sup |\nabla u|\ \1 _{B _{ \rho_u}}(x)\ \1 _{\{ |\xi|\le \delta\}}.
\end{align*}

All what follows relies on Lemma \ref{res-40} and we use its notation.
\\
By our hypotheses \eqref{eq-100b} and \eqref{eq-100c}, the constants $c''(\alpha),C''_ \rho( \alpha),c'_ \rho( \alpha), C'_ \rho( \alpha)$ attached to $ \Psi( \alpha)$ for all $ \alpha\in \mathcal{A},$ are uniform:
\begin{align*}
c'':= \sup _{t\in\ii, \alpha\in \mathcal{A}}c''(t,\alpha)< \infty,\qquad C''_ \rho:=\sup _{t\in\ii, x\in B _{ \rho+1},  \alpha\in \mathcal{A}}|\nabla_x ^2\Psi(t,x, \alpha)|< \infty,\\
c'_ \rho:=\inf _{t\in\ii, x\in B _{ \rho+1},  \alpha\in \mathcal{A}}|\nabla_x \Psi(t,x ,\alpha)|>0,\qquad 
C'_ \rho:=\sup _{t\in\ii, x\in B _{ \rho+1},  \alpha\in \mathcal{A}}|\nabla_x \Psi(t,x, \alpha)|< \infty.
\end{align*}

\begin{enumerate}[(i)]
\item
Control of $B_ \alpha$. Let us start estimating $\zeta$ defined at \eqref{eq-111}. We see with \eqref{eq-106} that $\nabla \hat\phi_\alpha (x)= \nabla^2 \psi_\alpha (x+z_x)$. With $C'' _{ \rho_u}< \infty$, this gives us
\begin{multline*}
\zeta(x) = |\det(\Id+\nabla\hat\phi_\alpha )|(x)-1 \\
\le  \Delta \psi_\alpha (x+z_x)
	+ \mathcal{R} _{ n-2}(\nabla^2 \psi_ \alpha(x+z_x))\  \sum _{ i,j} \big[\partial _{ ij}^2 \psi_ \alpha(x+z_x)\big]^2, 
\quad \forall x \in B _{ \rho_u},
\end{multline*}
where  $ \mathcal{R} _{ n-2}(\partial _{ ij}^2 \psi_ \alpha; 1\le i,j\le n) $ is a polynomial of order $n-2$. Consequently, with  \eqref{eq-105} and \eqref{eq-100db} we see that there exists some constant $c$ and a small $ \delta_*>0$ such that 
\begin{align*}
|\phi_\alpha (x)| \le \delta_*\implies \zeta(x)\le c |\phi_\alpha (x)|,\quad \forall  x\in B _{ \rho_u}.
\end{align*} 
Finally, we see that there exists  $C< \infty$ such that
\begin{align*}
B_ \alpha(x)\le C  \1 _{B _{ \rho_u}}(x)\ \1 _{\{ |\fa(x)|\le \delta_o(|x|)\}}\ 
k(x, \fa(x)) q(x)\ 
|\fa(x)|^2,
\quad \forall x \in\ZZ,\ \alpha\in \mathcal{A}.
\end{align*}
Note that $C$ and $ \delta_o$ do not depend on $ \alpha$ because all the bounds $c'', c'_ \rho,\dots$  are uniform in $ \alpha$.

\item
Control of $B^x_ \alpha$. We obtain a similar estimate: 
There exists  $C< \infty$ such that for all $ x \in\ZZ,$ $\alpha\in \mathcal{A}$
\begin{multline*}
B^x_ \alpha(x)\le C  \1 _{B _{ \rho_u}}(x)\ \1 _{\{ |\fa(x)|\le \delta_o(|x|)\}}
\big[|\nabla k_x(x,\fa(x))|q(x)+k(x,\fa(x))q(x)\\+k(x,\fa(x))|\nabla q(x)|\big]
\ |\fa(x)|^2.
\end{multline*}
We simply have to remark as above that $|\hat\phi_\alpha |\le 2|\phi_\alpha |$ on $B _{ \rho_u}$, whenever $ \delta_o>0$ is small enough.

\item
Control of $B^\xi_ \alpha$. We obtain a similar estimate: 
There exists  $C< \infty$ such that for all $ x \in\ZZ,$ $\alpha\in \mathcal{A}$
\begin{align*}
B^\xi_ \alpha(x)\le C  \1 _{B _{ \rho_u}}(x)\ \1 _{\{ |\fa(x)|\le \delta_o(|x|)\}}
\big[| \nabla_\xi k(x,\fa(x))| q(x)|\fa(x)|^3
		+k(x,\fa(x)) q(x) |\fa(x)|^2\big].
\end{align*}
This time the main estimate to invoke is \eqref{eq-103}.
\end{enumerate}
Finally the absolute value of the integrand in the right hand side of \eqref{eq-108} is upper bounded by
\begin{multline*}
C  \1 _{B _{ \rho_u}}(x)\ \1 _{\{ |\fa(x)|\le \delta_o(|x|)\}}\times\\
\Big(\big[\{k(x,\fa(x))+|\nabla k_x(x,\fa(x))|\}q(x)+k(x,\fa(x))|\nabla q(x)|\big]
\ |\fa(x)|^2\\
+  | \nabla_\xi k(x,\fa(x))| q(x)|\fa(x)|^3\Big)
\end{multline*}
for some constant $C$ 
and the identity \eqref{eq-108} gives for any $w_t\ge 0,$
\begin{align*}
 \lambda_t(w_t)
	&\le C \IZZs w(t,x) \1 _{B _{ \rho_u}}(x)\ \1 _{\{ |\xi|\le \delta_o(|x|)\}}
	\big[(k+|\nabla k|)(x,\xi) q(x)\ |\xi|^2\\
	&\hskip 4cm +k(x,\xi)|\nabla q(x)||\xi|^2
	+ | \nabla_\xi k(x,\xi)| q(x)|\xi|^3\big]\ \hat K _{ t,x}(d\xi)dx\\
	&\le CD\ \IZ w(t,x) \1 _{B _{ \rho_u}}(x)\,dx
\end{align*}
with 
\begin{multline*}
D:=\sup _{t\in\ii, x\in B _{ \rho_u}} \IZs \1 _{ \{|\xi|\le \delta_o(|x|)\}}
	\big[(k+|\nabla k|)(x,\xi) q(x)\ |\xi|^2
	+k(x,\xi)|\nabla q(x)||\xi|^2\\
	+ | \nabla_\xi k(x,\xi)| q(x)|\xi|^3\big]\ \hat K _{ t,x}(d\xi).
\end{multline*}
We see that 
\begin{align*}
\ell(w):=\Iii  \lambda_t(w_t)\, dt \le CD \IiZ w(t,x) \1 _{B _{ \rho_u}}(x)\,dtdx
\end{align*} is a bounded measure because it is assumed at \eqref{eq-119} and \eqref{eq-88}   that $D< \infty,$ and $C$ does not depend on $t$ because all the bounds $c'', c'_ \rho,\dots$  are uniform in $t$.
\end{proof}

\section{Entropic improvement}
\label{sec-entropy}

In this section, we extend the time reversal formula of Theorem \ref{res-36}  departing from the assumed regularity of $x\mapsto \Jf _{ t,x}(dy)$.

\subsection*{Entropic improvement} The strategy of this improvement is a variation on a theme  by Föllmer \cite{Foe85b,Foe86}. We present it in an abstract setting. It splits into two steps.

We start from some reference path measure $R\in\MO$  whose forward and backward generators $\LLf^R$ and $\LLb^R$ are known, for instance by means of Theorem \ref{res-36}. This strategy is aimed at deriving a time reversal formula  for any path \emph{probability} measure $P\in\PO$ with finite entropy 
\begin{align*}
H(P|R)< \infty
\end{align*}
with respect to $R$.

\subsubsection*{Step 1} \emph{Entropy, Girsanov and time reversal.}\ 
   Girsanov's theory tells us that under this finite entropy condition,
\begin{align*}
\LLf^P=\LLf^R+ \Gf ^{ P|R}
\end{align*}
and offers us an expression of $\Gf ^{ P|R}$. This  term is not regular in general: it is a measurable perturbation of $\LLf^R$ which satisfies some integrability condition (for the entropy to be finite). This really extends the class of path measure $P$ for which a time reversal formula can be derived.
\\
On the other hand, since the time reversal mapping is one-one,   the time reversal $P^*$ and $R^*$ of $P$ and $R$ satisfy 
\begin{align*}
H(P^*|R^*)=H(P|R)< \infty,
\end{align*}
so that Girsanov's theory applies as before, providing us with
$	%\begin{align*}
\LLf ^{ P^*}=\LLf ^{ R^*}+ \Gf ^{ P^*|R^*}.
$	%\end{align*}
 Since $\LLb^P_t=\LLf ^{ P^*} _{ 1-t}$, we arrive at
\begin{align*}
\LLb ^{ P}=\LLb ^{ R}+ \Gb ^{ P|R},
\end{align*}
with $\Gb ^{ P|R}_t:=\Gf  _{ 1-t}^{ P^*|R^*}.$  
Note that this shows that $P$ admits a backward generator, a property which is not granted in general, see \cite{W82, JP88}.

\subsubsection*{Step 2} \emph{IbP formula.}\ 
At this stage, we only know the existence of  $\Gb ^{ P|R}$ and its general shape. We still have to relate it precisely to $\Gf ^{ P|R}.$ The finiteness of $H(P^*|R^*)$ might imply (it will, under some additional hypotheses) that $$\IiZ |\LLb^P u|\,d\pb< \infty.$$ This is precisely  the specific assumption of part (a) of Theorem \ref{res-02} which gives us the IbP formula, i.e.\ an explicit expression of the backward generator in terms of the forward generator and the carré du champ.

\subsection*{Girsanov theory}

Let $R\in\MO$ be such that
\begin{align*}
\LLf^Ru=\Bf ^{ R, \delta}\scal\nabla u+\IZs [u(y)-u(x)-\nabla u(x)\scal \floo yx]\, \Jf _{ t,x}^R(dy),
\quad u\in C^2_c(\ZZ).
\end{align*}
In addition, it is assumed to be the unique solution to its martingale problem   in the sense of \eqref{eq-115}. 

For any $P\in\PO$ satisfying $H(P|R)< \infty,$ we have: $\CcZ\subset \dom\LLf^P,$ and  there exists a measurable function $j :\iZ\times \ZZ\to[0, \infty)$ such that for any $u\in\CcZ,$ and all $0\le t\le T, x\in\ZZ,$
\begin{align}\label{eq-30}
\LLf^P_tu(x)
	=\LLf^R_tu(x)+\IZs  [u(y)-u(x)](j  (t,x;y)-1)\Jf^R _{ t,x}(dy),
\end{align}
and 
\begin{align}\label{eq-28}
H(P|R)
	=H(P_0|R_0)
	+  \int _{ \ii\times\ZZ\times \ZZ} \hh ( j ) \, d[\pb\Jf^R]< \infty,
\end{align}
where $\pb\Jf^R(dtdxdy):=dt\pp_t(dx)\Jf^R _{ t,x}(dy),$ and $\hh$ is defined at \eqref{eq-56}.
For detailed proofs, see \cite{Leo11a}.
The function 
\begin{align}\label{eq-80}
\theta (a):= \hh(|a|+1)=(|a|+1)\log(|a|+1)-|a|,
\end{align}
 is a Young function satisfying  $ \theta (a)=a^2/2+ o _{ a\to 0}(a^2)$ and  $ \lim  _{ |a|\to \infty}\theta (a)/(|a|\log|a|)=1.$ Therefore, the Orlicz space $L^ {\theta }$ is equal to $L^2\cap L\log L .$ Since $ \theta (a)\le \hh(a+1),$ the estimate \eqref{eq-28} implies that
\begin{align*}
j -1\in L^ {\theta }(\pb\Jf^R).
\end{align*}
This justifies that  the integral in  \eqref{eq-30} is well-defined for any $u\in \CcZ .$ Otherwise stated, \eqref{eq-30} means that 
\begin{align}\label{eq-114}
\begin{split}
\Bf ^{ P, \delta}_t(x)&=\Bf ^{ R, \delta}_t(x) + \IZs  \floo yx\ (j  (t,x;y)-1)\Jf^R _{ t,x}(dy),\\
\Jf _{ t,x} ^P(dy)&=j(t,x;y)\, \Jf _{ t,x}^R(dy).
\end{split}
\end{align}
We see in particular that when the sample paths have bounded variations, taking $ \delta=0,$  we have $\Bf^P=\Bf^R$. Only the jump kernel is modified.

%\subsection*{dhflskdhgsmdfkhgsmlkdfghsmlkgj}\input{hyp}

\subsection*{Preliminary estimates}

From now on, we assume that  $R\in\MO$ satisfies the Hypotheses \ref{ass-08}.

So far, we used only Hypotheses \ref{ass-08}-a-b-c.
The supplementary Hypotheses \ref{ass-08}-d are designed for next lemmas to hold true.

\begin{lemma}[Control of $\Kb^R$]\label{res-39}
Let $R$ satisfy the Hypotheses \ref{ass-08}.\\
Then,
\begin{align}\label{eq-91}
 \sup _{t\in[t_o,T], x:|x|\le \rho}\IZs (1 \wedge |\xi|^2)\, \Kb _{ t,x}^R(d\xi)< \infty,
 \quad \forall  \rho\ge 0.
 \end{align}
\end{lemma}

\begin{proof}
By Theorem \ref{res-36}, $R$ satisfies the flux identity \eqref{eq-65}.

\subsubsection*{Under Hypothesis (d1)}
Since
\begin{align*}
\Jb^R _{ t,x}(dy)= \frac{d\Jf^R _{ t,y}}{d\rr_t} (x) \rr_t(dy)
	= \frac{d \tilde \pi_t}{d\pi_t}(x,y)\Jf^R _{ t,x}(dy)
	\overset{ \eqref{eq-81}}\le c( \rho)\Jf^R _{ t,x}(dy),
	\quad \forall x: |x|\le \rho,
\end{align*}
where $c( \rho):=\sup  \left\{
	 \frac{d\tilde \pi_t}{d\pi_t}(x,y); t\in\ii, x: |x|\le \rho, y\in\ZZ
	\right\} ,$
we obtain \eqref{eq-91} directly because
\begin{multline*}
\sup _{t\in\ii, x: |x|\le  \rho}\int    ( 1\wedge |\xi|)^2\, \Kb^R _{ t,x}(d\xi)
	\le c( \rho)\sup _{t\in\ii, x: |x|\le  \rho}\int   ( 1\wedge |\xi|)^2\, \Kf^R _{ t,x}(d\xi)\overset{ \eqref{eq-32}}< \infty.
\end{multline*}

\subsubsection*{Under Hypothesis (d2)}
By \eqref{eq-83} with $\sigma(x,y)=1 \wedge |y-x|^2,$ 
we have to verify
\begin{align}\label{eq-92}
 \sup _{t\in[t_o,T] y:|y|\le \rho}\IZs  \frac{d \rr^\xi_t}{d\rr_t}(y)\,(1\wedge|\xi|^2)\, \Kf ^R_{t, y-\xi}(d\xi)< \infty,
 \quad \forall  \rho\ge 0.
\end{align}
Our assumption \eqref{eq-89} is 
$	%\begin{align*}
\displaystyle{
c(\rho):=\sup _{t\in[t_o,T], y: |y|\le \rho, |\xi|\le \Delta( \rho)}
\frac{d \rr^\xi_t}{d\rr_t}(y)< \infty.
}
$	%\end{align*}
It implies
\begin{multline*}
 \sup _{t\in[t_o,T], y:|y|\le \rho}\IZs  \frac{d \rr^\xi_t}{d\rr_t}(y)\,(1\wedge|\xi|^2)\, \Kf ^R_{t, y-\xi}(d\xi)\\
 	\le c( \rho) \ \sup _{t\in[t_o,T], y:|y|\le \rho}\IZs (1\wedge|\xi|^2)\, \Kf ^R_{t, y-\xi}(d\xi)< \infty,
\end{multline*}
which is finite by hypothesis \eqref{eq-32}.
\end{proof}

Remark that for the proof of Lemma \ref{res-39}, it is enough in the unbounded variations case that  $r$ is  continuous, rather than $C ^{ 0,1} _{ t,x}$ as in \eqref{eq-88}.

\begin{lemma}[Control of $\Btb$] \label{res-42}
Assume that $R$ satisfies the Hypotheses \ref{ass-08}-a-b-c.
\\
Then,
\begin{align}\label{eq-95}
 \sup _{t\in[t_o,T], x:|x|\le \rho} |\Btb(t,x)|< \infty,
 \quad \forall  \rho\ge 0.
 \end{align}
\end{lemma}

\begin{proof}
If the sample paths have finite variation, the result is immediate because $\Btb =\Bb ^0=-\Bf$ and $\Bf$ is assumed to be locally bounded.
\\
Let us suppose now that \eqref{eq-97} holds. 
Since $\Btf$ is assumed to be locally bounded, and $r$ is assumed to be continuous and positive, it is sufficient to show that
\begin{align}\label{eq-96}
 \sup _{t\in[t_o,T], x:|x|\le \rho} r(t,x) |\Btf(t,x)+\Btb(t,x)|< \infty,
 \quad \forall  \rho\ge 0.
\end{align}
By \eqref{eq-83} and \eqref{eq-67},
\begin{multline*}
r(t,x)\,( \Btf+\Btb)(t,x)
	=r(t,x)\,\IZs \xd \Big( \Kf _{ t,x}(d\xi)- \frac{d \rr ^{ \xi}}{d\rr}(x)\, \Kf_{t, x-\xi}(d\xi)\Big)\\
	 = r(t,x)\,\IZs \xd \Big( \Kf _{ t,x}(d\xi)- \frac{r(x-\xi)}{r(t,x)}\, \Kf_{t, x-\xi}(d\xi)\Big)\\
	=\IZs \xd \Big(r(t,x) \Kf _{ t,x}(d\xi)- r(x-\xi)\Kf_{t, x-\xi}(d\xi)\Big)\\
	=\int _{ \mathcal{A}}\Big\{ r(t,x)\flo{\phi _{ t,x}( \alpha)}- r(x-\phi _{ t,x}( \alpha))\flo{\phi _{ t,x-\phi _{ t,x}( \alpha)}( \alpha)}\Big\}\, \Lambda( d \alpha);
\end{multline*}
Denoting $r(x)=r(t,x),$  $w=x-\phi _{ t,x}( \alpha),$  $\varphi(x)=\flo{\phi _{ t,x}( \alpha)}$, $\hat\varphi(x)=\flo{\hat\phi _{ t,x}( \alpha)}$  and $ \varphi(w)=\flo{\phi _{ t,x-\phi _{ t,x}( \alpha)}( \alpha)}$,  the integrand is
\begin{align*}
r(x) \varphi(x)-r(w) \varphi(w)
	&=r(x) \varphi(x)+r(w) \hat\varphi(x)\\
	&=r(x)[ \varphi(x)+\hat \varphi(x)]+[r(w)-r(x)] \hat\varphi(x)\\
	&=r(x)[ \varphi(x)+\hat \varphi(x)]+\nabla r(\tilde x)\scal (w-x)\, \hat\varphi(x)\\
	&=r(x)[ \varphi(x)+\hat \varphi(x)]+\nabla r(\tilde x)\scal \hat\varphi(x)\, \hat\varphi(x)
\end{align*}
for some  $\tilde x\in[w,x]$, where we used \eqref{eq-102}, \eqref{eq-101} and our hypothesis \eqref{eq-88}.
With \eqref{eq-105}  and \eqref{eq-103}, this gives
\begin{align*}
|r(x) \varphi(x)-r(w) \varphi(w)|
	\le c_ \rho\,| \varphi(x)|^2,\quad \forall x \in B_\rho,
\end{align*}
for some finite constant $c_\rho$ cooked up with $\sup _{ x\in B_ \rho}r(x),$ $ c'_\rho,$ $ C'_ \rho,$ and $C''_ \rho.$ Finally
\begin{align*}
r(x)\,|\Btf+\Btb|(t,x)
	\le c_ \rho \int _{ \mathcal{A}} \big| \flo{\phi _{ t,x}( \alpha)}\big|^2\ \Lambda(d \alpha)
	=c_ \rho \IZs |\xd|^2\, \Kf _{ t,x}(d \xi).
\end{align*}
We conclude with the assumption \eqref{eq-32}   that \eqref{eq-96} is true.
 \end{proof}

\subsection*{Applying the entropic improvement}

Next lemma asserts  that, when $H(P|R)< \infty,$ the backward jump kernel of $P$ is
$$
\Jb^P=j^* \Jb^R
$$
for some  function $j^*,$ and also that $P$ meets the hypotheses of Theorem \ref{res-02}-(a). This lemma completes the proof of Theorem \ref{res-43}.

If \eqref{eq-81} holds, i.e.\ under assumption (d1), we put   $t_o=0.$ 
Otherwise  under assumption (d2), $t_o$ is the same as  in  Hypothesis \eqref{eq-89}.

\begin{lemma}\label{res-18}
Assume that $R\in\MO$ satisfies the Hypotheses \ref{ass-08}, \eqref{eq-115}  and that in restriction to $[t_o,T]$ for some $0\le t_o<T,$ $P\in\PO$ has a finite entropy: $H(P _{ [t_o,T]}|R _{ [t_o,T]})< \infty.$ 
\begin{enumerate}[(a)]
\item
Then, $\CcZ\subset \dom \LLb^P$ and there exists some measurable function 
$j^* :[t_o,T]\times\ZZ\times\ZZ\to[0,\infty)$ satisfying  $ j^* -1\in L ^ {\theta }(\pb\Jb^R)$ and such that for any $u\in\CcZ$  
\begin{align}\label{eq-90}
\begin{split}
\LLb^P _tu(x)=\LLb^R_t u(x) +\IZ [u(y)-u(x)](j^* (t,x;y)-1)&\Jb^R _{ t,x}(dy),
\\ &t_o\le t\le T, x\in\ZZ,
\end{split}
\end{align}
and
\begin{align*}
H(P _{ [t_o,T]}| R _{ [t_o,T]})=H(P_T|R_T)+\int _{ [t_o,T]\times\ZZ} \hh(j^* )\, d [\pb\Jb^R]<\infty.
\end{align*}

\item
Moreover,    $\LLb^P u$ is integrable with respect to $ \1 _{ [t_o,T]}\pb.$
\end{enumerate}
\end{lemma}

Last estimate (b) is the main assumption of  Theorem \ref{res-02}-(a). The additional Hypotheses \ref{ass-08}-(d) on $R$ are introduced  to allow its proof.

\begin{proof}
\boulette{(a)}
As $H(P^*|R^*)=H(P|R)< \infty,$
the existence and integrability of $j^* $,   as well as the expression of the relative entropy follow from Girsanov theory.

\Boulette{(b)}
It remains to prove  the integrability of the backward generator $\LLb^Pu$.  Lemmas \ref{res-39} and \ref{res-42} show that $\LLb^R u$ is a bounded function on $ [t_o,T]\times\ZZ$, which is Hypothesis \ref{ass-08}-d3. Hence, we are left with controlling $$\IZ [u(y)-u(x)](j^* (t,x;y)-1)\,\Jb^R _{ t,x}(dy).$$
This is done by means of Young's inequality: $|ab|\le \theta(a) +\theta^*(b),$ where $ \theta$ is given at \eqref{eq-80} and $ \theta^*(b)=e ^{ |b|}-|b|-1$ its convex conjugate.   For any $ x,y$ 
$
|u(y)-u(x)|
	\le c_u  D_u(x,y)
$
where $c_u:=\max(2 \sup|u|,\sup|\nabla u|)$ and $D_u(x,y):= \1 _{ \supp u}(x) ( 1\wedge | y-x|).$ By Young's inequality,
\begin{multline*}
\Big|\int _{ [t_o,T]\times\ZZZ} [u(y)-u(x) ](j^* (t,x;y)-1)\pp_t(dx)\Jb^R _{ t,x}(dy)dt\Big|\\
\le c_u \int D_u\ |j^*-1|\, \, d[\pb\Jb^R]\\
\le c_u \Big\{ \int  \theta^*(D_u)\, d[\pb\Jb^R]+ \int  \theta(j^*-1)\, d[\pb\Jb^R]\Big\}.
\end{multline*}
The last integral is finite because $\int  \theta(j^*-1)\, d[\pb\Jb^R]\le 3\int\hh(j^*)\,d[\pb\Jb^R]< \infty$, as a consequence of the finite entropy assumption $H(P _{[ t_o,T]}|R _{ [t_o,T]})< \infty$. 
\\
Let us prove that the remaining integral  is also finite. We have 
\begin{multline*}
 \int\theta^*(D_ u)\, d[\pb\Jb^R]
 	\le T \sup _{t,x}\IZ \theta^*(D_ u)\, d\Jb _{ t,x}^R\\
	\le cT  \sup _{(t,x)\in[t_o,T]\times \supp u}\IZ (1 \wedge |y-x|^2)\, \Jb _{ t,x}^R(dy)
\end{multline*}
 where $c>0$ is a large enough constant such that: $ \theta^*(1 \wedge d)\le c (1 \wedge d^2),$ for all $d\ge 0,$ (note that $ \theta^*(b)=b^2/2+ o_{ b\to 0}(b^2)).$  We conclude with  \eqref{eq-91}.
\end{proof}

%\bibliographystyle{plain}
%\bibliography{bib-christian.bib}

\end{document}